\documentclass[a4paper,12pt]{article}
\usepackage{amssymb}
\usepackage{amsmath}
\usepackage{amsthm}
\usepackage{mathrsfs}
\usepackage[frame,curve,knot]{xy}
\xyoption{all}
\def\dj{d\kern-.30em\raise1.25ex\vbox{\hrule width .3em height .03em}}
\def\Dj{\rlap{\kern-.70em\raise0.75ex\vbox{\hrule width .3em height .03em}}} 

\newtheorem{theorem}{Theorem}[section]
\newtheorem{proposition}[theorem]{Proposition}
\newtheorem{lemma}[theorem]{Lemma}
\newtheorem{cor}[theorem]{Corollary}

\newtheorem{definition}[theorem]{Definition}

\newtheorem{remark}[theorem]{Remark}

\newcommand{\ad }{\textrm{ad}}
\newcommand{\ABM }{{}^\cA_\B \hspace{-.1cm}\mathcal{M}}
\newcommand{\ABMB }{{}^\cA_\B \hspace{-.1cm}\mathcal{M}_\B}
\newcommand{\ann}{\mbox{Ann}}
\newcommand{\AqM }{{}^{\Aqr}\hspace{-.1cm}\mathcal{M}}
\newcommand{\AqMB }{{}^{\Aqr}\hspace{-.1cm}\mathcal{M}_\B}

\newcommand{\Aqr }{\overleftarrow{\cA }}
\newcommand{\B }{\mathcal{B}}

\newcommand{\cA }{\mathcal{A}}
\newcommand{\cB }{\mathcal{B}}
\newcommand{\cF}{\mathcal{F}}

\newcommand{\cO}{\mathcal{O}}

\newcommand{\cqg}{\mathbb{C}_q[G]}
\newcommand{\cqgl}{\mathbb{C}_q[G/L_S]}

\newcommand{\C }{\mathbb{C}}
\newcommand{\Cu}{\underline{C}}
\newcommand{\del }{\partial}
\newcommand{\delb }{{\overline{\partial}}}
\newcommand{\dif }{\mathrm{d}}

\newcommand{\fun}{\underline{\,\,\,}}
\newcommand{\Gam }{\varGamma }

\newcommand{\Gduw}[1][]{\Gamma _{\del ,\mathrm{u}}^{\wedge #1}}
\newcommand{\Gdbuw}[1][]{\Gamma _{\delb ,\mathrm{u}}^{\wedge #1}}
\newcommand{\Gdifuw}[1][]{\Gamma _{\dif ,\mathrm{u}}^{\wedge #1}}

\newcommand{\gfrak}{\mathfrak{g}}

\newcommand{\gr}{\mathrm{gr}}
\newcommand{\gu}{\underline{g}}
\newcommand{\hfrak}{\mathfrak{h}}
\newcommand{\hb}{\overline{h}}
\newcommand{\Hom}{\mathrm{Hom}}
\newcommand{\hght}{\mathrm{ht}}
\newcommand{\hu}{\underline{h}}
\newcommand{\id}{\mathrm{Id}}

\newcommand{\im}{\mathrm{Im}}
\newcommand{\invol}{\eta}

\newcommand{\kopr }{\varDelta }
\newcommand{\kow }{\varDelta }

\newcommand{\lfrak}{\mathfrak{l}}

\newcommand{\Lin }{\mathrm{Lin}}

\newcommand{\MK}{{\mathcal M}_{\uqls}}

\newcommand{\N }{\mathbb{N}}
\newcommand{\nfrak}{\mathfrak{n}}
\newcommand{\ob}{\overline{\phantom{a}}}

\newcommand{\op }{\mathrm{op}}
\newcommand{\oep}{\overline{\vep}}
\newcommand{\ophi}{\overline{\varphi}}
\newcommand{\Om }{\Omega }
\newcommand{\omvol }{\omega_{vol}}
\newcommand{\ot }{\otimes }
\newcommand{\pair }[2]{\langle #1,#2 \rangle }

\newcommand{\pfrak}{\mathfrak{p}}

\newcommand{\psp}{P_S^+}

\newcommand{\ract }{\triangleleft}

\newcommand{\roots }{R }

\newcommand{\slfrak}{\mathfrak{sl}}

\newcommand{\sroots }{\pi }
\newcommand{\thetab}{\overline{\theta}}

\newcommand{\Tor }{\mbox{Tor}}

\newcommand{\U }{U}
\newcommand{\ufrak }{\mathfrak{u}}
\newcommand{\ug }{U(\mathfrak{g})}

\newcommand{\uqg }{U_q(\mathfrak{g})}
\newcommand{\uqgc }{\check{U}_q(\mathfrak{g})}

\newcommand{\uql }{U_q(\mathfrak{l})}
\newcommand{\uqls }{U_q(\mathfrak{l}_S)}
\newcommand{\uqlsm }{U_q(\mathfrak{l}_S^-)}
\newcommand{\uqlsp }{U_q(\mathfrak{l}_S^+)}
\newcommand{\uqps }{U_q(\mathfrak{p}_S)}

\newcommand{\uqpsop }{U_q(\mathfrak{p}_{S}^{\mathrm{op}})}

\newcommand{\uqnp }{U_q(\mathfrak{n}^+)}
\newcommand{\uqnm }{U_q(\mathfrak{n}^-)}
\newcommand{\uqslz }{U_q(\mathfrak{sl}_2)}
\newcommand{\ubar }{\overline{\U }}
\newcommand{\wght }{\mathrm{wt}}
\newcommand{\wlat }{P}
\newcommand{\wurz }{\pi }
\newcommand{\Vcal }{\mathcal{V} }
\newcommand{\vep }{\varepsilon }

\newcommand{\WM }{W^{M(\lambda)} }

\newcommand{\vphiu}{\underline{\varphi }}

\newcommand{\Wcal }{\mathcal{W} }

\newcommand{\Z}{\mathbb{Z}}

\begin{document}

\title{Differential forms via the Bernstein-Gelfand-Gelfand resolution
for quantized irreducible flag manifolds}

\author{Istv\'an Heckenberger 
        and Stefan Kolb\footnote{supported by the German Research Foundation (DFG)}}

\date{November 16, 2006}

\maketitle

\begin{abstract}
 The quantum group version of the Bernstein-Gelfand-Gelfand resolution
 is used to construct a double complex of $\uqg$-modules with exact
 rows and columns. The locally finite dual of its total complex is
 identified with the de Rham complex for quantized irreducible flag manifolds.    
\end{abstract}


\section{Introduction}
Over the last two decades a vast amount of papers has been devoted to
the translation of classical geometric concepts to
coordinate algebras appearing in the theory of quantum groups. It is a
recurring theme that such constructions are possible if the underlying
geometric object can be expressed in purely Lie algebraic terms. 
A list of examples where this translation has a very simple and compelling form might
include the standard definition of the $q$-deformed coordinate algebra
$\cqg$ inside the dual Hopf algebra of $\uqg$ \cite[9.1.1]{b-Joseph} or
the construction of the quantum group version of the homogeneous
coordinate ring of a flag
manifold \cite[9.1.6]{b-Joseph}. Certainly, one always aims for quantum
effects, as for instance Drinfeld duality, which transcend the
classical undeformed situation. However, we will not encounter
significant quantum effects in this paper.

Differential forms are an example of a geometric concept where the translation from the
classical to the quantum group setting is far from obvious in
general. However, there is a notion of covariant differential calculus
on quantum spaces, introduced by S.L.~Woronowicz \cite{a-Woro2}, which has
attracted much attention for many years (\cite{b-KS} and references
therein). It soon turned out that for a general quantum spaces there
exists no canonical construction of a covariant differential calculus.
However, in \cite{a-heko06} we showed that for quantized irreducible flag
manifolds $G/P_S$ where $G$ is a simple complex affine algebraic group
and $P_S$ a standard parabolic subgroup there exists a $q$-analog of
the de Rham complex which in many respects behaves like its undeformed
counterpart. The aim of the present paper is to relate this complex
to its Lie algebraic shadow, the Bernstein-Gelfand-Gelfand (BGG)
resolution. In the quantum case such a construction was suggested by
L.L.~Vaksman and a first indication of its feasibility can be found in
\cite{a-SinVa} where generalized Verma modules are used to obtain
$q$-analogs of differential one forms. 

The main result of the present paper, Theorem \ref{mainResult}, states that
the de Rham complex investigated in \cite{a-heko06} can also be
obtained as the locally finite dual of a BGG-like sequence of
$\uqg$-modules induced by $\uqls$-modules, where $\lfrak_S$ denotes the
Levi factor of the parabolic subalgebra $\pfrak_S\subset \gfrak$.
More precisely the BGG resolution for quantum groups \cite{a-heko06p} is used to define
quantum analogs of the complexes of holomorphic and antiholomorphic
differential forms on flag manifolds (Proposition \ref{DelCalcIso} and
Section \ref{delbcalc}). In Section \ref{difcalc} we introduce a double
complex the rows and columns of which are closely related to the BGG
resolutions used to obtain the holomorhpic and antiholomorphic
differentials, respectively. The desired de Rham complex is then
obtained as the locally finite dual of the total complex of this
double complex.

The reason why we have to consider $\uqls$-modules, instead of
$\uqps$-modules as one might expect, lies in the definition of the
coordinate algebra $\cqgl$ describing the quantum flag manifold. Its
classical counterpart is the coordinate ring of the affine algebraic
variety $G/L_S$ where $L_S$ denotes the Levi factor of $P_S$. The
advantage of this approach lies in the fact that $\cqg$ is a
Hopf-Galois extension of $\cqgl$. Thus M.~Takeuchi's categorical
equivalence \cite{a-Tak79} applies and one can make use of results on
differential calculi on quantum homogeneous spaces \cite{a-HK-QHS}. 

A result similar in spirit has recently been obtained in
\cite{a-SSV06p}. In that paper the universal higher order differential
calculus constructed in \cite{a-SinVa} is identified with the category $\cO$
dual of the $q$-version of the BGG-resolution. Hence in the approach
taken in \cite{a-SSV06p} Takeuchi's categorical equivalence is not
available and the authors have to revert to specialization techniques.
In the present paper, on the other hand, all results are proved for
any deformation parameter $q\in \C$ which is not a root of unity.

As the reader might at first be put off by the technical nature of our
paper we now state the main result in the special case of one dimensional
quantum complex projective space also known as standard
quantum sphere. This simplest example of an irreducible quantized flag
manifold in itself has been subject to
various publications, e.g.~\cite{a-DabSit03}, \cite{a-SchmWag04},
\cite{a-Maj05}. We believe that our analysis will lead to new insight even
in this simplest case. 

Recall that $\uqslz$ denotes a Hopf algebra generated by elements $E$,
$F$, $K$, and $K^{-1}$ and relations given for instance in
\cite[3.1]{b-KS}. Let $\uql$ denote the subalgebra generated by $K$
and $K^{-1}$, for $n\in \Z$ let $V(n)$ denote the one-dimensional $\uql$-module
generated by one element $v_n$ with the action $Kv_n=q^{2n} v_n$, and
define $W(n,m):=\uqslz\ot_{\uql} V(n{-}m)$. Note that $W(0,0)$ is a coalgebra
and that $W(n,m)$ is a left and right comodule over $W(0,0)$ with
coactions given by 
\begin{align*}
  u\ot v_n\mapsto (u_{(1)}{\ot} v_0)\ot (u_{(2)}{\ot} v_n) \quad \mbox{and}\quad
  u\ot v_n\mapsto (u_{(1)}{\ot} v_n)\ot (u_{(2)}{\ot} v_0),
\end{align*}
respectively, where Sweedler notation is used. Consider the following
sequence of $\uqslz$-modules, $W(0,0)$-bicomodules
 \begin{align} \label{Wsequence}
   \xymatrix{& & W(1,0)\ar[rd]^{\varphi_{1,0;0}}& &\\
              0\ar[r]&W(1,1)\ar[ru]^{\varphi_{1;1,0}}\ar[rd]^{{\varphi}_{1,0;1}}
                     &\bigoplus &W(0,0)\ar[r]&0  \\
              & & W(0,1)\ar[ru]^{\varphi_{0;0,1}}&&
              }
  \end{align}
where $\varphi_{a;b,c}(u\ot v_{a-b})=uE\ot v_{a-c}$ and
$\varphi_{a,b;c}(u\ot v_{a-c})=uF\ot v_{b-c}$. 
The locally finite dual of $W(a,b)$ is defined by
\begin{align*}
  \Om^{a,b}=\{f\in W(a,b)^\ast\,|\,\dim(f\uqslz)<\infty\}
\end{align*}
where $W(a,b)^\ast$ denotes the linear dual space of $W(a,b)$.
As $W(0,0)$ is a $\uqslz$-module coalgebra the space $\cB=\Om^{0,0}$ is a
$\uqslz$-module algebra and as such $\cB$ coincides with the standard
quantum sphere. Dualizing (\ref{Wsequence}) one obtains a sequence of
$\uqslz$-module $\cB$-bimodules
\begin{align*}
  \xymatrix{& & \Om^{1,0}\ar[rd]^{\del_{1;1,0}}& &\\
              0\ar[r]&\cB\ar[ru]^{\del_{1,0;0}}\ar[rd]_{{\del}_{0;0,1}}
                     &\bigoplus &\Om^{1,1}\ar[r]&0. \\
              & & \Om^{0,1}\ar[ru]_{\del_{1,0;1}}&&
              }
\end{align*}
The main result of this paper, Theorem \ref{mainResult}, states in this
special case, that this sequence coincides with the well known de Rham
complex \cite{a-Po92} over the standard quantum sphere $\cB$.
As an application one can for instance immediately read off the twisted
cyclic cocycle calculated in \cite[Lemma 4.4]{a-SchmWag04}.

We now describe the contents of each section of this paper in some
detail. In Section \ref{prel} we fix notations. Moreover, we compare
the standard resolution of the trivial module with the parabolic
version of the BGG resolution and show that these two coincide if
$\gfrak/\pfrak_S$ is irreducible. This result should be well known but
we were not able to track it in the literature. On the one hand it
explains once again why it is necessary to assume irreducibility of
the considered flag manifolds. On the other hand, it implies that
certain weights $w.0$ are incomparable in the Bruhat order. 

Section \ref{QG} serves purely to fix notations for quantum groups and
to recall M.S.~K\'eb\'e's results on triangular decompositions of $\uqg$
with respect to parabolic subalgebras. In Section \ref{QVerma} we
quickly review the $q$-analog of the BGG resolution which by
\cite{a-heko06p} is exact if $q$ is not a root of unity. Section
\ref{induced} is devoted to $\uqg$-modules induced by irreducible
$\uqls$-modules. We denote the category of finite direct sums of such
modules by $\Wcal$. In Subsection \ref{Homest} we derive technical
properties of standard maps between objects in $\Wcal$ related to the
BGG resolution. The locally finite duals of objects in $\Wcal$ are
interpreted as yet another realization of Takeuchi's categorical
equivalence in Section \ref{CatEquiv}.

The main technical work is done in the final Section \ref{diffdual}.
First the main results from \cite{a-heko06} are recalled. Then the
differential calculi $(\Gamma_\del,\del)$, $(\Gamma_\delb, \delb)$, and
$(\Gamma_\dif, \dif)$ from that paper are interpreted as
locally finite duals of BGG-like sequences in $\Wcal$. 

Explicit calculations flooded by symbols are inherent to proofs in 
quantum group theory. For the convenience of the reader we have collected all
commonly used notation in order of appearance in an appendix.

This project started out during a two week visit of the first author at Virginia
Tech in March 2005. He wishes to thank the mathematics department of
Virginia Tech for hospitality. The second author was a
DFG postdoctoral fellow at Virginia Tech from August 2004 to July
2006. He is very grateful to the mathematics department, and to his
hosting researcher Gail Letzter in particular, for encouragement
and support. Both authors are indepted to Leonid Vaksman for leading them towards
quantized irreducible flag manifolds and the BGG resolution. 

\section{Preliminaries}\label{prel}
Let $\N$, $\Z$, and $\C$ denote the positive integers, the integers,
and the complex numbers, respectively. We write $\N_0$ to denote the
nonnegative integers.
\subsection{Notations}\label{notation1}
First, to fix notations some general notions related to Lie algebras are
recalled. Let $\gfrak$ be a finite dimensional complex simple Lie algebra of rank $r$ and
let $\hfrak\subset \gfrak$ be a fixed Cartan subalgebra. Let $R\subset \hfrak^\ast$ 
denote the root system associated with $(\gfrak,\hfrak)$.
Choose an ordered basis $\wurz=\{\alpha_1,\dots,\alpha_r\}$ of simple roots
for $R$ and let $R^+$ (resp.~$R^-$) be the set of positive (resp.~negative)
roots with respect to $\wurz$.
Moreover, let $\gfrak=\nfrak^+\oplus \hfrak\oplus \nfrak^-$ be the
corresponding triangular decomposition. 
Identify $\hfrak$ with its dual via the Killing form. The induced
nondegenerate symmetric bilinear form on $\hfrak^*$
is denoted by $(\cdot,\cdot)$. The root lattice $Q=\Z R$
is contained in the weight lattice $P=\{\lambda\in\hfrak^\ast\,|\,
(\lambda,\alpha_i^\vee)\in\Z\,\forall \alpha_i\in\wurz\}$ where 
$\alpha_i^\vee:=2\alpha_i/(\alpha_i,\alpha_i)$. In order to avoid roots of the deformation
parameter $q$ in the following sections we rescale $(\cdot,\cdot)$ such that
$(\cdot,\cdot):P\times P\rightarrow \Z$.
For $\mu,\nu\in \wlat$ we write $\mu\ge \nu$ if $\mu-\nu$ is a sum of positive roots.
The height function $\hght:Q\rightarrow \Z$ is defined by 
$\hght(\sum_{i=1}^r n_i\alpha_i)=\sum_{i=1}^r n_i$.

The fundamental weights $\omega_i\in\hfrak^\ast$, $i=1,\dots,r$, are
characterized by $(\omega_i,\alpha_j^\vee)=\delta_{ij}$.
Let $P^+$ denote the set of dominant integral weights, i.~e.~the $\N_0$-span of
$\{\omega_i\,|\,i=1,\dots,r\}$.
Recall that $(a_{ij}):=(2(\alpha_i,\alpha_j)/(\alpha_i,\alpha_i))$ is the
Cartan matrix of $\gfrak$ with respect to $\wurz$. We will 
write $Q^+=\N_0R^+$.

For $\mu\in P^+$ let $V(\mu)$ denote the finite dimensional 
irreducible $\gfrak$-module of highest weight $\mu$.
Moreover, let $\Pi (V(\mu))$ denote the set of weights of the 
$\gfrak$-module $V(\mu)$.

Let $G$ denote the connected simply connected complex Lie group with Lie
algebra $\gfrak$. For any set $S\subset \wurz$ of simple roots define 
$Q_S=\Z S$, $Q_S^+=Q_S\cap Q^+$, and $R_S^\pm:=Q_S\cap R^\pm$.
Let $P_S$ and $P_S^{\mathrm{op}}$ denote the corresponding standard parabolic
subgroups of $G$ with Lie algebra
\begin{align}
  \pfrak_S=\hfrak\oplus\bigoplus_{\alpha\in R^+\cup R^-_S}\gfrak_\alpha,\qquad
  \pfrak^{\mathrm{op}}_S=\hfrak\oplus\bigoplus_{\alpha\in R^-\cup R_S^+}
  \gfrak_\alpha,
\end{align}
respectively. Moreover,
\begin{align*}
  \lfrak_S:=\hfrak\oplus\bigoplus_{\alpha\in R^+_S\cup R^-_S}
\gfrak_\alpha
\end{align*}
is the Levi factor of $\pfrak_S$ and $L_S=P_S\cap P_S^\op\subset G$ denotes
the corresponding subgroup. 

The generalized flag manifold $G/P_S$ is called irreducible if the
adjoint representation of $\pfrak_S$ on $\gfrak/\pfrak_S$ is irreducible.
Equivalently, $S=\pi\setminus\{\alpha_i\}$ where $\alpha_i$ appears in any
positive root with coefficient at most one. For a complete list of all
irreducible flag manifolds consult e.g.~\cite[p.~27]{b-BastonEastwood}.
Note that the irreducible flag manifolds coincide with the irreducible compact
Hermitian symmetric spaces \cite[X\S 6.3]{b-Helga78}

Define
$\psp :=\{\lambda\in P \,|\,(\lambda,\alpha_i)/d_i \in
  \N_0\,\forall\alpha_i \in S\}$.
To $\lambda\in \psp $ we associate the finite dimensional, irreducible
$\lfrak_S$-module $M(\lambda)$ of highest weight $\lambda$.

Let $W$ denote the Weyl group of $\gfrak$ generated by the
reflections corresponding to the simple roots in $\pi$. 
For any $\alpha\in R^+$  let $s_\alpha\in W$ denote the reflection
on the hyperplane orthogonal to $\alpha$ with respect to $(\cdot,\cdot)$. 
Let $W_S\subset W$ denote the subgroup generated by the 
reflections corresponding to simple roots in $S$. 
Moreover, define
\begin{align*}
  W^S=\{w\in W\,|\, R_S^+\subset w R^+\}.
\end{align*}
By a well known result of B.~Kostant any element $w\in W$ can
be decomposed uniquely in the form $w=w_Sw^S$ where $w_S\in W_S$ and
$w^S\in W^S$. Moreover, if $l$ denotes the length function on $W$ then
this decomposition satisfies $l(w)=l(w_S)+l(w^S)$.

The following technical Lemma will be used in the proof of Propositions
\ref{Xi-prop} and \ref{TakLokFin}.

\begin{lemma}\label{KinU}
For any $\lambda\in \psp\cap \wlat$ there exist $\mu\in \wlat^+$
which allows an injective $\lfrak_S$-module map 
$M(\lambda)\hookrightarrow V(\mu)$.
\end{lemma}

\noindent {\bf Proof:} Choose $w\in W$ such that 
$\mu:=w^{-1}\lambda\in \wlat^+$.
Write $ w=w_S\,w^S$ where $w_S\in W_S$ and $w^S\in W^S$. Then
$w^S\mu=w_S^{-1}\lambda$. Let $v_{w^S\mu}\in V(\mu)$ denote a nonzero
vector of weight $w^S\mu$. Note that $v_{w^S\mu}$ is a highest weight
vector for $\lfrak_S$. Indeed, if $w^S\mu+\alpha_i\in \Pi(V(\mu))$ 
for some $\alpha_i\in S$ then 
$(w^S)^{-1}(w^S\mu+\alpha_i)=\mu+(w^S)^{-1}\alpha_i\notin 
\Pi(V(\mu))$ since $(w^S)^{-1}\alpha_i\in \roots^+$.
Therefore $w_S^{-1}\lambda\in \psp$ and $\lambda\in \psp$ and
hence $w^S\mu=w_S^{-1}\lambda=\lambda$.\,
$\blacksquare$

\vspace{.5cm}

Recall that the shifted action of the Weyl group $W$ on $P$ is defined
in terms of the ordinary Weyl group action by
\begin{align*}
  w.\mu=w(\mu+\rho)-\rho
\end{align*}
where $\rho$ is half the sum of all positive roots or equivalently
$\rho=\sum_{i=1}^r \omega_i$.
Moreover, for $w,w'\in W$ write $w\rightarrow w'$ if there exists
$\alpha\in R^+$ such that $w=s_\alpha w'$ and $l(w)=l(w')+1$. The Bruhat
order $\le $ on $W$ is then given by the relation
\begin{align*}
w\le w' \Leftrightarrow \
&\text{there exists } n\ge 1,\ w_2,\ldots ,w_{n-1}\in W,\\
&\text{such that }w=w_1\rightarrow w_2\rightarrow \ldots
\rightarrow w_n=w'.
\end{align*}

\subsection{Standard-resolution and BGG-resolution} 
Let $\gfrak$ be a complex Lie algebra and $\pfrak$ a subalgebra. 
In \cite{a-BGG75} 
I.~N.~Bernstein, I.~M.~Gelfand, and S.~I.~Gelfand
have given the following generalization of the standard resolution of
Lie algebra cohomology. The adjoint action of $\pfrak$ on 
$\gfrak/\pfrak$ endows each exterior product 
$\Lambda^k(\gfrak/\pfrak)$ with the structure of a $U(\pfrak)$-module.
Define 
\begin{align*}
  D_k=\U(\gfrak)\ot_{U(\pfrak)}\Lambda^k(\gfrak/\pfrak)
\end{align*}
and 
\begin{align*}
  d_0: U(\gfrak)\ot_{U(\pfrak)}\C=D_0\rightarrow \C,\quad 
   u\ot x\mapsto \vep(u)x
\end{align*}
where $\vep$ denotes the counit of $U(\gfrak)$. 
Moreover, for $k\ge 1$ define operators $d_k: D_k\rightarrow D_{k-1}$
in the following way.
Let $X_1,\dots, X_k$ be elements of $\gfrak/\pfrak$. Let
$Y_1,\dots, Y_k\in \gfrak$ be arbitrary representatives of 
$X_1,\dots, X_k$, respectively, and put
\begin{align}
  d_k(&X\ot X_1\wedge\dots\wedge X_k)= 
  \sum_{i=1}^k(-1)^{i+1} (XY_i\ot X_1\wedge\dots\wedge \hat{X}_i\wedge
  \dots\wedge X_k)\label{d_kdef}\\
  &+\sum_{1\le i<j\le k}(-1)^{i-j}(X\ot\overline{[Y_i,Y_j]}\wedge
  X_1\wedge\dots\wedge\hat{X_i}\wedge\dots\wedge\hat{X_j}\wedge\dots
  \wedge X_k).\nonumber
\end{align}
Here $X\in U(\gfrak)$ and we write $\overline{Y}$ for the image of
the element $Y\in\gfrak$ in $\gfrak/\pfrak$. Moreover, $\hat X$ denotes omission of the
element $X$.
One obtains a complex
\begin{align*}
  D_\ast: 0\leftarrow \C \stackrel{d_0}{\leftarrow}D_0
  \stackrel{d_1}{\leftarrow}D_1\stackrel{d_2}{\leftarrow}D_2
  \stackrel{d_3}{\leftarrow}\, \dots
\end{align*}
which is exact by \cite[Thm.~9.1]{a-BGG75}. 
In general the complex $D_\ast$ does not have an analogue for
quantum universal enveloping algebras.

Let now $\gfrak$ be a finite dimensional simple complex Lie algebra 
and  $\pfrak_S\subset\gfrak$ a standard parabolic subalgebra
as in the previous subsection. 

For any irreducible highest weight module $V(\mu)$ of $\gfrak$, 
where $\mu\in P^+$, in generalization of 
\cite{a-BGG75} J.~Lepowsky \cite{a-Lep77}
constructed an exact sequence of $U(\gfrak)$-modules
\begin{align*}
  0\leftarrow V(\mu) \stackrel{}{\leftarrow}C_0
  \stackrel{}{\leftarrow}C_1\leftarrow\,\dots\leftarrow 
  C_{\dim(\gfrak/\pfrak_S)}\leftarrow 0
\end{align*}
where
\begin{align*}
  C_n=\bigoplus_{\makebox[0cm]{$w\in W^S,\atop l(w)=n$}}
         U(\gfrak)\ot_{U(\pfrak_S)}M(w.\mu).
\end{align*}
Here the differentials are given as linear combinations of standard
maps of the occurring generalized Verma modules.
In particular if $\mu=0$ one obtains an exact sequence
\begin{align}\label{BBG0}
  C_\ast: 0\leftarrow \C \stackrel{\delta_0}{\leftarrow}C_0
  \stackrel{\delta_1}{\leftarrow}C_1
  \stackrel{\delta_2}{\leftarrow}\,\dots\leftarrow 
  C_{\dim(\gfrak/\pfrak_S)}\leftarrow 0.
\end{align}
For general parabolics the sequences of $U(\gfrak)$-modules
$C_\ast$ and $D_\ast$ are not isomorphic. Indeed, if $\gfrak/\pfrak_S$ is not irreducible
then not even $D_1$ and $C_1$ need to be isomorphic. However, one has
the following result. 
\begin{proposition}\label{C=D}
  If $\gfrak/\pfrak_S$ is irreducible then the corresponding 
  complexes of 
  $U(\gfrak)$-modules $(C_\ast,\delta_\ast)$ and $(D_\ast,d_\ast)$ 
  are isomorphic.
\end{proposition}

\noindent
{\bf Proof:} Consider the Lie subalgebra 
\begin{align*}
  \ufrak_S^-=\bigoplus_{\alpha\in\roots^-\setminus \roots_S^-}
  \gfrak_{\alpha}\subset \gfrak.
\end{align*}
One has decompositions $\gfrak=\ufrak_S^-\oplus \pfrak_S$ and
$\ug\cong \U(\ufrak_S^-)\ot U(\pfrak_S)$ by the Poincar\'e-Birkhoff-Witt
Theorem. Note that both $D_\ast$ and $C_\ast$ are free resolutions
of the trivial left $U(\ufrak_S^-)$-module $\C$. Thus both sequences can
be used to compute $\Tor^{U(\ufrak_S^-)}_j(\C,\C)$ where the first
entry $\C$ denotes the trivial right $U(\ufrak_S^-)$-module.
Note that if $\gfrak/\pfrak_S$ is irreducible then the second
term in (\ref{d_kdef}) vanishes, because $\ufrak_S^-$ is commutative. 
Thus in the complex 
$\C\ot_{U(\ufrak_S^-)}D_\ast$ all differentials vanish and therefore
\begin{align*}
  \dim(\Tor^{U(\ufrak_S^-)}_j(\C,\C))=\dim(\Lambda^j(\gfrak/\pfrak_S)).
\end{align*}
Similarly the sequence (\ref{BBG0}) yields 
(cp.~\cite[Cor.~3.11]{a-Lep77})
\begin{align*}
  \dim(\Tor^{U(\ufrak_S^-)}_j(\C,\C))=
  \sum_{\makebox[0cm]{$w\in W^S,\atop l(w)=j$}}\dim (M(w.0)).
\end{align*}
Thus one obtains
\begin{align}\label{dimequal}
  \dim(\Lambda^j(\gfrak/\pfrak_S))=
  \sum_{\makebox[0cm]{$w\in W^S,\atop l(w)=j$}}\dim (M(w.0)).
\end{align}
As $\Lambda^j(\gfrak/\pfrak_S)$ and $\ug$ are graded by the root lattice
and the differentials $d_j$ of the complex $D_\ast$ respect this
grading one can define a $\Z$-grading of the complex $D_\ast$ by
\begin{align*}
  \deg(u\ot v)=(\omega_s,\wght(u)+\wght(v))=(\omega_s,\wght(u))-j
\end{align*}
where $u\in \ug$ and $v\in \Lambda^j(\gfrak/\pfrak_S)$ are homogeneous
elements. Similarly the complex $C_\ast$ is $\Z$-graded by the
same formula where now $v\in M(w.0)$ for some $w\in W^S$ with $l(w)=j$.
Assume now that the complexes $D_\ast$ and $C_\ast$ are isomorphic
as $\Z$-graded complexes of $\ug$-modules up to complex-degree
$k$. This holds for $k=0$. Set $Z_k:=\ker d_k\subset C_k=D_k$.
Then $d_{k+1}(D_{k+1})=Z_k=\delta_{k+1}(C_{k+1})$ by exactness of 
the sequences. Moreover, from (\ref{d_kdef}) one obtains that
$d_{k+1}$ is injective when restricted to 
$1\ot \Lambda^{k+1}(\gfrak/\pfrak_S)$.
Note that 
\begin{align*}
  1\ot \Lambda^{k+1}(\gfrak/\pfrak_S)=\{x\in D_{k+1}\,|\,
  \deg(x)=-(k+1)\}
\end{align*}
and hence
\begin{align*}
  d_{k+1}(1\ot \Lambda^{k+1}(\gfrak/\pfrak_S))=Z_k^{-(k+1)}
\end{align*}
for $Z_k^{-(k+1)}:=\{x\in Z_k\,|\,\deg(x)=-(k+1)\}$. 
As $Z_k$ does not contain any element $x$ such that $\deg(x)>-(k+1)$
one has $\delta_{k+1}(\ug\ot_{U(\pfrak_S)}M(w.0))=0$ if 
$(\omega_s,\wght(w.0))>-(k+1)$. By the injectivity of $d_{k+1}$ and
(\ref{dimequal}) one has  
\begin{align*}
  \dim(Z_k^{-(k+1)})=\sum_{\makebox[0cm]{$w\in W^S,\atop l(w)=k+1$}}
  \dim(M(w.0)).
\end{align*}
As $\delta_{k+1}$ maps onto $Z_k$ this implies
\begin{align*}
  \delta_{k+1}
  \Big(1\ot \bigoplus_{\makebox[0cm]{$w\in W^S,\atop l(w)=k+1$}}
  M(w.0)\Big) =Z_k^{-(k+1)}.
\end{align*}
Hence one obtains in view of (\ref{dimequal}) that 
the composition
\begin{align*}
  \phi_{k+1}: 
  1\ot \bigoplus_{\makebox[0cm]{$w\in W^S,\atop l(w)=k+1$}}
  M(w.0) 
  \stackrel{\delta_{k+1}}{\longrightarrow} Z^{-(k+1)}_k
  \stackrel{(d_{k+1})^{-1}}{\longrightarrow}
  1\ot \Lambda^{k+1}(\gfrak/\pfrak_S)
\end{align*}
is an isomorphism of $U(\pfrak_S)$-modules. By construction
\begin{align*}
  \id\ot \phi_{k+1}: \ug\ot _{U(\pfrak_S)}\bigoplus_{\makebox[0cm]{$w\in W^S,\atop l(w)=k+1$}}
  M(w.0) \rightarrow \ug\ot _{U(\pfrak_S)}\Lambda^{k+1}(\gfrak/\pfrak_S)
\end{align*}
extends the isomorphism of complexes to degree $k+1$.
$\blacksquare$

\vspace{.5cm}

As an application one obtains the following Corollary.
\begin{cor}\label{uncomparable}
  Let $\gfrak/\pfrak_S$ be irreducible and let
  $w_1,w_2\in W^S$ be elements of equal length $l(w_1)=l(w_2)$.
  
  1) One has $w_1.0-w_2.0 \in Q_S$.
  
  2) If $w_1\neq w_2$ then $w_1.0-w_2.0\notin Q_S^+$.
  
  \noindent Moreover, if $w,w'\in W^S$ and $l(w)=l(w')+1$ then
  $\omega_s(w.0-w'.0)=1$.
\end{cor}
\noindent
{\bf Proof:}
 1) By the above Proposition $w_1.0$ and $w_2.0$ occur as weights of
  $\Lambda^{l(w_1)}\gfrak/\pfrak_S$. As $\gfrak/\pfrak_S$ is
  irreducible the weights of $\gfrak/\pfrak_S$ differ by
  elements in $Q_S$. Then so do the weights of
  $\Lambda^{l(w_1)}\gfrak/\pfrak_S$.
  
 2) Assume $w_1.0-w_2.0\in Q_S^+$, or equivalently
$w_1\rho-w_2 \rho\in Q_S^+$. Multiplication by $w_2^{-1}$
and the definition of $W^S$ yield
\begin{align*}
  w_2^{-1}w_1\rho-\rho\in Q^+.
\end{align*}
Since $\rho$ is dominant and $W$ acts faithfully on $\rho$ 
one obtains a contradiction unless $w_1=w_2$.

The last statement follows from the fact that the map $\phi_k$ from the
proof of the proposition is an isomorphism of $U(\pfrak_S)$-modules.
$\blacksquare$

\medskip

 \begin{remark}\upshape 1) In the above corollary the condition of irreducibility of 
$\gfrak/\pfrak_S$ can't be dropped. Indeed, for $S=\emptyset$ one has $Q_S=\{0\}$ but 
$w_1.0\neq w_2.0$ for $w_1\neq w_2$.

2) Also one can't replace $0$ by a more general weight $\mu\in P^+$. Consider for
example $\gfrak=\mathfrak{sl}_4$, $S=\{\alpha_1,\alpha_3\}$, and $\mu=\omega_3$.
Then $s_2s_3.\omega_3=\omega_3-2\alpha_3-3\alpha_2$ and 
$s_2s_1.\omega_3=\omega_3-\alpha_1-2\alpha_2$ and hence
$s_2s_3.\omega_3-s_2s_1.\omega_3=\alpha_1-\alpha_2-2\alpha_3\notin Q_S$.
On the other hand $s_2s_3$ and $s_2s_1$ are elements of $W^S$ of equal length.
\end{remark}

\section{Quantum groups}\label{QG}

\subsection{Definition of $\uqg$ and $\cqg$}\label{uqgcqg}
We keep the notations of the previous section. 
Let $q\in \C\setminus\{0\}$ be not a
root of unity. The $q$-deformed universal enveloping algebra $\uqg$ 
associated to $\gfrak$ is considered here as the complex algebra 
generated by elements $K_i,K_i^{-1},E_i,F_i$, $i=1,\dots,r$, 
and relations as given for instance in \cite[6.1.2]{b-KS}.
In particular one has
\begin{align}\label{relations}
\begin{aligned}
&\begin{aligned}
 K_iE_j&=q^{(\alpha_i,\alpha_j)}E_jK_i,& K_iF_j&=q^{-(\alpha_i,\alpha_j)}
 F_jK_i,
\end{aligned}&\\
&\begin{aligned}
 E_iF_j-F_jE_i&=\delta_{ij}\frac{K_i-K_i^{-1}}{q^{(\alpha_i,\alpha_i)/2}-q^{-(\alpha_i,\alpha_i)/2}},
\end{aligned}&
\end{aligned}
\end{align}
The algebra $\uqg$ has a Hopf algebra structure with coproduct given by
\begin{align*}
  \kopr K_i&= K_i\otimes K_i,&\kopr E_i&= E_i\otimes K_i+1\otimes E_i,&
  \kopr F_i&= F_i\otimes 1 + K_i^{-1}\otimes F_i\nonumber.
\nonumber
\end{align*}
These formulae for the coproduct imply in particular that the antipode
$\kappa$ of $\uqg$ is given by
\begin{align*}
  \kappa(K_i)=K_i^{-1},\qquad\kappa(E_i)=-E_iK_i^{-1},
  \qquad\kappa(F_i)=-K_iF_i.
\end{align*}
The counit will be denoted by $\vep$. We will make frequent use of Sweedler notation
in the form $\kow u=u_{(1)}\ot u_{(2)}$ for $u\in \uqg$. Moreover, for any $u,x\in \uqg$ we 
will write $(\ad u)x=u_{(1)}x\kappa(u_{(2)})$ to denote the left adjoint action.

There exists a uniquely determined algebra isomorphism coalgebra antiisomorphism $\invol$ 
of $\uqg$ such that
\begin{align*}
\invol(E_i)=F_i,\quad \invol(F_i)=E_i,\quad 
\invol(K_i^{\pm 1})=K_i^{\mp 1}. 
\end{align*}

Let $\uqnp,\uqnm\subset \uqg$ denote the subalgebras generated
by $\{E_i\,|\,1\le i \le r\}$ and $\{F_i\,|\,1\le i \le r\}$, 
respectively. Let $U^0\subset \uqg$ be the subalgebra generated by
$\{K_i, K_i^{-1}\,|\,1\le i \le r\}$.
Moreover, let $G_+\subset \uqg$ denote the subalgebra
generated by $\{E_iK_i^{-1}\,|\,1\le i \le r\}$.

For $\mu\in P^+$  let $V(\mu)$ denote the uniquely determined finite
dimensional irreducible left $\uqg$-module with highest weight $\mu$.
More explicitly, there exists a highest weight vector
$v_\mu\in V(\mu)\setminus\{0\}$ satisfying
\begin{align}
  E_iv_\mu=0,\quad
  K_iv_\mu=q^{(\mu,\alpha_i)}v_\mu
  \qquad \mbox{ for all } i=1,\dots,r.
\end{align}
In general a vector $v\in V(\mu)$ is called a weight vector of weight
$\wght (v)\in P$ if $K_i v=q^{(\wght (v),\alpha_i)}v$ 
for all $i=1,\dots,r$.

The dual $V^\ast$ of a finite dimensional $\uqg$-module $V$ is 
defined as the dual vector space with the $\uqg$-action given by
\begin{align*}
  (uf)(v)=f(\kappa(u)v)\qquad \forall\,v\in V,\,f\in V^\ast,\,u\in \uqg.
\end{align*}

For any left $\uqg$-modules $V$ define a new $\uqg$-module $V_\eta$ to be the same 
vector space with the left $\uqg$-module structure $\bullet_\eta$ given by
\begin{align}\label{involact}
  u\bullet_\invol v:=\invol(u)v\qquad\mbox{for all } u\in\uqg, 
  v\in V.
\end{align}
Note that $V(\mu)_\eta\cong V(\mu)^\ast$. 

As usual the $q$-deformed coordinate ring $\cqg$ is defined to be the 
subspace of the linear dual $\uqg^\ast$ spanned by the matrix 
coefficients of the finite dimensional irreducible representations 
$V(\mu)$, $\mu\in P^+$.
For $v\in V(\mu)$, $f\in V(\mu)^\ast$ the matrix coefficient
$c^\mu_{f,v}\in \uqg^\ast$ is defined by
\begin{align*}
  c^\mu_{f,v}(X)=f(Xv).
\end{align*}
The linear span of matrix coefficients of $V(\mu)$
\begin{align}
  C^{V(\mu)}=\Lin_\C\{c^\mu_{f,v}\,|\,v\in V(\mu), f\in V(\mu)^\ast\}
\end{align}
obtains a $\uqg$-bimodule structure by
\begin{align}
  (Yc^\mu_{f,v}Z)(X)=f(ZXYv)=c^\mu_{fZ,Yv}(X).
\end{align}
Here $V(\mu)^\ast$ is considered as a right $\uqg$-module. Note that by
construction
\begin{align}\label{peterweyl}
  \cqg\cong\bigoplus_{\mu\in P^+} C^{V(\mu)}
\end{align}
is a Hopf algebra and the pairing
\begin{align}\label{pair}
  \cqg\ot \uqg\rightarrow \C
\end{align}
is nondegenerate.

\subsection{Nilpotent and parabolic subalgebras}\label{kebe}
For $S\subset \sroots$ let $\uqls\subset \uqg$ denote the Hopf 
subalgebra generated by $E_i,F_i, K_j, K_j^{-1}$ for all 
$\alpha_i\in S$ and all $j$.
Moreover, let $V_-\subset \uqg$ denote the subalgebra generated by
the elements of the set
\begin{align*}
  \{(\ad k)F_i\,|\, k\in \uqls, \alpha_i\notin S  \}.
\end{align*}
Analogously, let $V_+\subset \uqg$ denote the subalgebra generated by
the elements of the set
\begin{align*}
  \{(\ad k)(E_iK_i^{-1})\,|\, k\in \uqls, \alpha_i\notin S  \}.
\end{align*}  
As $(\ad E_i) F_j=0=(\ad F_i)(E_jK_j^{-1})$ for all $i\neq j$
one has $V_-\subset \uqnm$ and $V_+\subset G_+$. 
By \cite[Prop.~4.2]{a-Kebe99} multiplication gives isomorphism
\begin{align*}
  V_-\ot \uqlsm\rightarrow \uqnm,\qquad V_+\ot \uqlsp\rightarrow G_+
\end{align*}
where $\uqlsm:=\uqnm\cap \uqls$ and $\uqlsp:=G_+\cap \uqls$. Thus the 
triangular decomposition $\uqg\cong \uqnm\ot G_+\ot U^0$ yields
\begin{align}
  \uqg&\cong V_-\ot \uqlsm\ot V_+\ot \uqlsp \ot U^0\nonumber\\
      &\cong V_-\ot V_+\ot \uqlsm\ot \uqlsp\ot U^0\nonumber\\
      &\cong V_-\ot V_+\ot \uqls.\label{V-V+uqls}
\end{align}
Here in the second line the isomorphism
\begin{align*}
  \uqlsm\ot V_+ \rightarrow V_+\ot \uqlsm, \qquad k\ot v\mapsto 
  (\ad\, k_{(1)})v\ot k_{(2)}
\end{align*}
is used, and the last line uses the triangular decomposition of 
$\uqls$.
In a similar manner one obtains
\begin{align}
    \uqg\cong V_+\ot V_-\ot \uqls.\label{V+V-uqls}
\end{align}

The parabolic subalgebra $\uqps\subset \uqg$ is defined by
\begin{align*}
  U_q(\pfrak_S):=\langle E_i,K_i, F_j\,|\, \alpha_i\in \pi,
                            \alpha_j\in S \rangle.
\end{align*}  
Note that $\uqps$ coincides with the subalgebra generated by
$\uqls$ and $V_+$. Thus by (\ref{V-V+uqls}) multiplication
yields isomorphisms
\begin{align}
  V_+\ot \uqls &\cong \uqps,\label{v+ls}\\
  V_-\ot \uqps &\cong \uqg.\label{v-ps}
\end{align}

\section{Quantum generalized Verma modules}\label{QVerma}

\subsection{Notation}\label{QVnots}
For $\lambda\in \psp $ as in the classical case $q=1$ let $M(\lambda)$
denote the finite dimensional, irreducible
$\uqls$-module of highest weight $\lambda$.
Note that $M(\lambda)$ can be turned into an $\uqps$-module by
setting $E_i v=0$ for all generators $E_i$, $\alpha_i\notin S$, and
$v\in M(\lambda)$.

\begin{definition}\label{VMlamdef}
  For $\lambda\in \psp $, define the quantum generalized Verma module
  $V^{M(\lambda)}$ by
  \begin{align*}
    V^{M(\lambda)}:=\uqg\ot_{\uqps} M(\lambda). 
  \end{align*}  
  If $S=\emptyset$ and $\lambda\in P$ we will write
  $V^\lambda:=V^{M(\lambda)}.$
\end{definition}  
Note that by (\ref{v-ps}) one has isomorphisms of $U^0$-modules 
$V^\lambda\cong \uqnm \ot \C^\lambda$ and
$V^{M(\lambda)}\cong V_-\ot M(\lambda)$ where $\C^\lambda$ denotes the
one-dimensional $U^0$-module of weight $\lambda$. 

Note moreover that $V^{M(\lambda)}_\eta\cong V^{M(\lambda)^\ast}$.
Indeed, let $\xi_{-\lambda}\in M(\lambda)^\ast$ denote the up to
scalar multiplication uniquely determined element of weight 
$-\lambda$. Then $1\ot \xi_{-\lambda}\in V^{M(\lambda)^\ast}$ is
a cyclic vector and a set of relations determining $V^{M(\lambda)^\ast}$ is given by
\begin{align*}
  K_j^{\pm 1}(1\ot \xi_{-\lambda})=q^{\mp (\lambda,\alpha_j)}
  1\ot \xi_{-\lambda},\quad 
  E_i^{(\lambda,\alpha_i^\vee)+1}(1\ot \xi_{-\lambda})=0
\end{align*}
for all $\alpha_i\in S$ and for all $j$.
The same relations hold for the cyclic vector 
$1\ot v_\lambda\in V_\invol^{M(\lambda)}$.

\vspace{.5cm}

\noindent{\bf Remark:} The notations used here slightly differ from 
the original notations in \cite{a-Lep77}. Recall that $\rho$ denotes
half the sum of the positive roots and define
\begin{align*}
  \rho_S:=\frac{1}{2}\sum_{\alpha\in R^+_S} \alpha.  
\end{align*}
Note that for all $\lambda\in P$ one has
\begin{align*}
  \lambda\in \psp  \Leftrightarrow \lambda-\rho+\rho_S\in \psp .
\end{align*}
J.~Lepowsky considered modules
obtained by twisted induction in the classical case $q=1$ and defined
$V^{M(\lambda)}:=U(\gfrak)\ot_{U(\pfrak_S)} M(\lambda-\rho+\rho_S)$.
Translation between the two settings is straightforward.

\vspace{.5cm}
 
Let $\lambda\in \psp $ and let $v_\lambda\in M(\lambda)$ 
denote a vector of weight $\lambda$. 
For any $\uqg$-module homomorphism
$g:V^{M(\lambda)}\rightarrow V^{M(\mu)}$ there exists an element
$F\in \uqnm$ such that $g(1\ot v_{\lambda})=F\ot v_{\mu}$.
Note that $F$ is uniquely determined up to addition of an element in the 
annihilator of $1\ot v_{\mu}\in V^{M(\mu)}$.
We will say that the homomorphism $g$ \textit{is determined by} $F$.

\subsection{The Bernstein-Gelfand-Gelfand resolution}\label{BGGsection}
We now briefly recall the quantum analog of the Bernstein-Gelfand-Gelfand resolution.
This construction has been in detail considered in \cite{a-heko06p} for $q$ not a root of 
unity.

Fix a dominant integral weight $\mu\in \wlat^+$. For all 
$j=0,\dots, \dim(\gfrak/\pfrak_S)$ define
\begin{align*}
  C_j^S:=\bigoplus_{w\in W^S, l(w)=j} V^{M(w.\mu)}.
\end{align*}
Note that $V^{M(w.\mu)}$ is a highest
weight module with highest weight $w.\mu $. Therefore
$V^{M(w.\mu)}$ is a natural quotient of $V^{w.\mu }$.

As in \cite[Section 4]{a-Lep77} one constructs $\uqg$-module maps 
$\varphi^S_j:C_j^S\rightarrow C_{j-1}^S$
for all $j=1, \dots, \dim(\gfrak/\pfrak_S)$.
More explicitly, for all $w\in W$, fix an embedding 
$V^{w.\mu}\subset V^\mu$. Then for all $w,w'\in W$ with $w\le w'$
one has a fixed embedding $f_{w,w'}:V^{w.\mu}\rightarrow V^{w'.\mu}$.

A quadruple $(w_1,w_2,w_3,w_4)$ of elements of $W$ is called a square
if $w_2\neq w_3$ and 
\begin{align*}
  w_1\rightarrow w_2\rightarrow w_4, \qquad
  w_1\rightarrow w_3\rightarrow w_4.
\end{align*}
By \cite[Lemma 10.4]{a-BGG75} to each arrow 
$w_1\rightarrow w_2$ ($w_1,w_2\in W$) one can 
assign a number $s(w_1,w_2)=\pm 1$ such that for every square, the
product of the numbers assigned to the four arrows occurring in it is 
$-1$.
Let $w,w'\in W^S$ such that $l(w)=l(w')+1$. If $w\rightarrow w'$ then
let $h_{w,w'}:V^{M(w.\mu)}\rightarrow V^{M(w'.\mu)}$ denote the
(standard) map induced by the map
\begin{align*}
  s(w,w')f_{w,w'}:V^{w.\mu}\rightarrow V^{w'.\mu}.
\end{align*}
Otherwise, define $h_{w,w'}=0$.
The differential $\varphi^S_j$ is now defined as the sum of all
$h_{w,w'}$ where $l(w)=j=l(w')+1$.

Note moreover, that for $\mu\in \wlat^+$ there exists a surjective map
of $\uqg$-modules
\begin{align}\label{epsdef}
  \vep_\mu:C_0^S=\uqg\ot_{\uqps} M(\mu)\rightarrow V(\mu),\qquad
  u\ot v_{\mu,M}\mapsto uv_{\mu,V}
\end{align}
where $v_{\mu,M}\in M(\mu)$ and $v_{\mu,V}\in V(\mu)$ denote vectors 
of weight $\mu$.

\begin{theorem}{\upshape \cite[Section 3.4]{a-heko06p}}
  The sequence
  \begin{align}\label{BGGres}
    0\longrightarrow C_{\dim \gfrak/\pfrak_S}^S
    \stackrel{\varphi^S_{\dim(\gfrak/\pfrak_S)}}{\longrightarrow} \cdots 
    \stackrel{\varphi^S_1}{\longrightarrow} C_0^S
    \stackrel{\vep_\mu}{\longrightarrow} V(\mu)\longrightarrow 0
  \end{align}
  is exact and $\varphi^S_j(V^{M(w.\mu)})\neq 0$
  for all $j=1,\dots,\dim(\gfrak/\pfrak_S)$ and all $w\in W^S$ with 
  $l(w)=j$.
\end{theorem}

\begin{remark}\label{standnonzero}
  \upshape In the quantum case the fact that for $w\rightarrow w'$ the standard
  map $h_{w,w'}:V^{M(w.\mu)}\rightarrow V^{M(w'.\mu)}$ is nonzero has
  not been explicitly stated in \cite{a-heko06p}. However, this
  property can be verified analogously to formula (1) in the proof of
  \cite[Lemma 9.2.14]{b-Kumar02}. The necessary fact that for
  $\mu,\lambda\in \hfrak^\ast$ the simple module $V(\mu)$ is a
  subquotient of $V^{\lambda}$ if and only if
  $\Hom(V^\mu,V^\lambda)\neq 0$ follows as in \cite{a-Neidhardt84}
  after translation of \cite[Sections 1-6]{a-RCW82} to the quantum case.
\end{remark}

By construction there exists $y_{w,w'}^\mu\in\uqnm$ such that 
$f_{w,w'}(u\ot v_{w.\mu})=u y^\mu_{w,w'}\ot v_{w'.\mu}$.
Thus in terms of the elements $y_{w,w'}^\mu\in\uqnm$ the map 
$h_{w,w'}$ is given by 
\begin{align*}
  h_{w,w'}(u\ot v_{w.\mu})=s(w,w')u y^\mu_{w,w'}\ot v_{w'.\mu}.
\end{align*}
In later considerations the main focus will be on the case 
$\mu=0$. In this case define $y_{w,w'}:=s(w,w')y^0_{w,w'}\in\uqnm$.

\section{$\uqg$-modules induced by $\uqls$-modules}\label{induced}
\subsection{Notation}\label{inducedNots}
\begin{definition}\label{WMlamdef}
  For $\lambda\in \psp $, define
  \begin{align*}
    W^{M(\lambda)}:=\uqg\ot_{\uqls} M(\lambda). 
  \end{align*}  
  If $S=\emptyset$ and $\lambda\in P$ we will write
$W^\lambda:=W^{M(\lambda)}.$
\end{definition}
Note that multiplication yields isomorphisms $W^\lambda\cong \uqnm\ot\uqnp\ot \C^\lambda$ 
and $W^{M(\lambda)}\cong V_-\ot V_+\ot M(\lambda)$
of $U^0$-modules.
Note moreover that $W^{M(\lambda)}$ is in general not a highest 
weight module.
In analogy to the observation after Definition \ref{VMlamdef}
one obtains
$W_\invol^{M(\lambda)}\cong W^{M(\lambda)^\ast}$.

\subsection{The functor $\fun:\Vcal\rightarrow \Wcal$}\label{functor}

By Definitions \ref{VMlamdef} and \ref{WMlamdef}
there exists a natural surjective $\uqg$-module homomorphism
\begin{align*}
  \Phi_\lambda: W^{M(\lambda)}\rightarrow V^{M(\lambda)}.
\end{align*}
\begin{proposition}\label{liftprop}
  For any $\uqg$-module homomorphism
  $g:V^{M(\lambda)}\rightarrow V^{M(\mu)}$ there exists a uniquely determined 
  $\uqg$-module homomorphism $\gu:W^{M(\lambda)}\rightarrow W^{M(\mu)}$ such that
  the diagram
  \begin{align}\label{diagram}
    \xymatrix{W^{M(\lambda)}\ar[r]^\gu \ar[d]^{\Phi_\lambda} &W^{M(\mu)}\ar[d]^{\Phi_\mu}\\
              V^{M(\lambda)}\ar[r]^g &V^{M(\mu)}
              }
  \end{align}
  commutes. 
\end{proposition}  
\noindent {\bf Proof:} Assume that $g$ is determined by $F\in \uqnm$ as in the end of 
Section~\ref{QVnots}. To obtain commutativity of the diagram (\ref{diagram}) one has to 
define $\gu(1\ot v_\lambda)=F\ot v_\mu$.
We have to check that $\gu$ is well defined. To this end consider
$0=u\ot v_\lambda\in W^{M(\lambda)}$, or 
equivalently $u\in \uqg \ann_{\uqls}(v_\lambda)$. We have to show that
$uF\in \uqg\ann_{\uqls} (v_\mu)$. 

Using the decomposition (\ref{V+V-uqls}) and the fact that $V_-\uqls F\subset V_-\uqls$ one 
may assume that $u\in V_-\uqls$.  The relation $g(\Phi_\lambda(u\ot v_\lambda))=0$ implies
$uF\in \uqg\ann_{\uqps}(v_\mu)$. Hence 
\begin{align*}
  uF&\in \uqg\ann_{\uqps}(v_\mu)\cap V_-\uqls \\ &\stackrel{(\ref{v-ps})}{=}
  V_-(\ann_{\uqps}(v_\mu)\cap\uqls)= V_-\ann_{\uqls}(v_\mu). \qquad \blacksquare
\end{align*}


Let $\Vcal$ and $\Wcal$ denote the full subcategory of the category of $\uqg$-modules
whose objects are finite direct sums of $\uqg$-modules $V^{M(\lambda)}$ and
$W^{M(\lambda)}$, where $\lambda\in P^+_S$, respectively.  By Proposition
\ref{liftprop} there exists a canonical functor
$\fun:\Vcal\rightarrow \Wcal$ such that
\begin{align*}
  \underline{\bigoplus_{i=1}^n V^{M(\lambda_i)}}= \bigoplus_{i=1}^n W^{M(\lambda_i)}.
\end{align*}

\begin{proposition}\label{PresProp}
  The functor $\fun:\Vcal\rightarrow \Wcal$ is exact.
\end{proposition}
\noindent{\bf Proof:} Recall that $V^{M(\lambda)}\cong V_-\ot M(\lambda)$
and $W^{M(\lambda)}\cong V_+\ot V_-\ot M(\lambda)$.
With respect to these decompositions one gets for any $V_1,V_2\in \Vcal$ and any 
$g:V_1\rightarrow V_2$ the relation $\gu=\id_{V_+}\ot g$.
Hence $\fun$ preserves exactness.
$\blacksquare$

\vspace{.5cm}

Let $\uqpsop\subset \uqg$ denote the subalgebra generated by
the elements $E_j,K_i,F_i$ for $\alpha_i\in \wurz$, $\alpha_j\in S$.
For any $\mu\in \wlat^+$ define a map
\begin{align*}
  \underline{\vep}_\mu:W^{M(\mu)}=\uqg\ot_{\uqls} M(\mu)
                                &\rightarrow \uqg\ot_{\uqpsop}V(\mu)\\
     \quad u\ot v_{\mu,M}&\mapsto u\ot v_{\mu,V}                           
\end{align*}
where as in (\ref{epsdef}) the symbols $v_{\mu,M}\in M(\mu)$ and 
$v_{\mu,V}\in V(\mu)$ denote vectors of weight $\mu$.
Then by the same argument as in the proof of Proposition 
\ref{PresProp} the BGG resolution
(\ref{BGGres}) induces an exact sequence
 \begin{align}\label{CundSeq}
    0\longrightarrow \Cu_{\dim \gfrak/\pfrak_S}^S
    \stackrel{\vphiu^S_{\dim(\gfrak/\pfrak_S)}}{\longrightarrow} 
    \cdots 
    \stackrel{\vphiu^S_1}{\longrightarrow} \Cu_0^S
    \stackrel{\underline{\vep}_\mu}{\longrightarrow} 
     \uqg\ot_{\uqpsop}V(\mu)\longrightarrow 0.
  \end{align}

\subsection{Homomorphisms and estimates}\label{Homest}
Note that for $\mu,\nu\in \psp$ the left $\uqls$-module $M(\mu)\ot M(\nu)^\ast$ is generated
by one element $v_\mu\ot \xi_{-\nu}$ where $v_\mu\in M(\mu)$ and $\xi_{-\nu}\in M(\nu)^\ast$
denote a highest and a lowest weight vector, respectively.
A complete set of relations for $M(\mu)\ot M(\nu)^\ast$ is given by
\begin{align}
  E_i^{(\nu,\alpha_i^\vee)+1}(v_\mu\ot \xi_{-\nu})&=0\nonumber\\
  F_i^{(\mu,\alpha_i^\vee)+1}(v_\mu\ot \xi_{-\nu})&=0\label{modul-rels}\\
  (K_j-q^{(\mu-\nu,\alpha_j)})(v_\mu\ot \xi_{-\nu})&=0\nonumber
\end{align}
where $\alpha_i\in S$ and $\alpha_j\in\pi$. This follows for instance
from \cite[Prop.~5.2]{a-Joseph99} using the fact that the module
generated by one element and the relations (\ref{modul-rels}) is
integrable and the generator is a cyclic weight vector. In Section \ref{difcalc} we will be
interested in homomorphisms between $\uqg$-modules induced by $\uqls$-modules  
$M(\mu)\ot M(\nu)^\ast$. Here we derive well definedness and some properties of such maps.

For $w,w'\in W^S$, $w\rightarrow w'$, and $\mu\in \wlat^+$ recall the definition of the 
element $y^\mu_{w,w'}\in \uqnm$ from Section \ref{BGGsection} and define
$x^\mu_{w,w'}:=\eta(y^\mu_{w,w'})$. Define $\uqls$-module homomorphisms
\begin{align*}
  \theta_2:M(w.\mu)&\rightarrow V^{M(w'.\mu)}\cong V_+V_-\ot M(w'.\mu),\\
  uv_{w.\mu}&\mapsto u y^\mu_{w,w'} \ot v_{w'.\mu},\\
  \thetab_2:M(w.\mu)^\ast&\rightarrow V^{M(w'.\mu)^\ast}\cong V_-V_+\ot M(w'.\mu)^\ast,\\
  u\xi_{-w.\mu}&\mapsto u x^\mu_{w,w'} \ot \xi_{-w'.\mu}.
\end{align*}
\begin{proposition}\label{triangProp}
  Let $w,w'\in W^S$, $w\rightarrow w'$, $\mu\in \wlat^+$, and $\nu\in \psp$. There 
  are uniquely determined injective $\uqls$-module homomorphisms
  \begin{align*}
    \theta_1:M(w.\mu)\ot M(\nu)^\ast &\rightarrow \uqg\ot_{\uqls}(M(w'.\mu)\ot M(\nu)^\ast),\\
    \thetab_1:M(\nu)\ot M(w.\mu)^\ast &\rightarrow \uqg\ot_{\uqls}(M(\nu)\ot M(w'.\mu)^\ast)
  \end{align*}
  such that
  \begin{align}
    \theta_1(v_{w.\mu}\ot \xi_{-\nu})&=y^\mu_{w,w'} \ot (v_{w'.\mu}\ot \xi_{-\nu}),\label{phi}\\
    \thetab_1(v_{\nu}\ot \xi_{-w.\mu})&=x^\mu_{w,w'} \ot (v_{\nu}\ot \xi_{-w'.\mu}).\label{phib}
  \end{align}
  Moreover, in $V_+V_-\ot M(w'.\mu)\ot M(\nu)^\ast$ one has for all weight vectors
  $v\in M(w.\mu)$, $\xi\in M(\nu)^\ast$
  \begin{align}\label{phi1(votxi)}
    \theta_1(v\ot \xi)\in \theta_2(v)\ot \xi +\sum_{\zeta<\wght(\xi)} V_-\ot M(w'.\mu)\ot 
      M(\nu)^\ast_\zeta.
  \end{align}
 Similarly, in $V_-V_+\ot M(\nu)\ot M(w'.\mu)^\ast$ one has for all weight vectors 
 $v\in M(\nu)$, $\xi\in M(w.\mu)^\ast$ 
  \begin{align}\label{phib1(votxi)}
    \thetab_1(v\ot \xi)\in P_{23}(\thetab_2(\xi)\ot v) 
    +\sum_{\zeta>\wght(v)} V_+\ot M(\nu)_\zeta \ot M(w'.\mu)^\ast
  \end{align} 
  where $P_{23}$ denotes the flip of the second 
  and the third tensor factor.
\end{proposition}

\noindent{\bf Proof:} 
The maps $\theta_1$ and $\thetab_1$ are uniquely determined by formulae (\ref{phi}) 
and (\ref{phib}), respectively. It remains to verify that they are well defined and 
injective. Fix $\alpha_i\in S$ and let $U_i\subset \uqg$ denote the subalgebra isomorphic
to $U_q(\slfrak_2)$ generated by $E_i$, $F_i$, and $K_i^{\pm 1}$. Note that
\begin{align}\label{F-rels}
  F_i^{(w.\mu,\alpha_i^\vee)+1} y^\mu_{w,w'}\ot (v_{w'.\mu}\ot \xi_{-\nu})=0
\end{align}
for all $\alpha_i\in S$. Indeed, as the standard map $h_{w,w'}$ is well defined one obtains
\begin{align*}
  F_i^{(w.\mu,\alpha_i^\vee)+1} y^\mu_{w,w'}\in (\uqg \ann_{\uqls} v_{w'.\mu})\cap \uqnm.
\end{align*}  
Hence (\ref{F-rels}) follows from the fact that $\xi_{-\nu}$ is a lowest weight vector.

Note that $\uqg\ot_{\uqls} (M(w.\mu)\ot M(\nu)^\ast)$ is an integrable $\uqls$-module. 
Hence the weight vector
$y^\mu_{w,w'}\ot (v_{w'.\mu}\ot \xi_{-\nu})$ can be written as a sum of weight vectors of
weight $w.\mu-\nu$ which generate pairwise nonisomorphic irreducible $U_i$-modules. 
By (\ref{F-rels}) among these irreducible $U_i$-modules there is one of lowest weight
$w.\mu-\nu-(w.\mu,\alpha_i^\vee)\alpha_i$ and all other $U_i$-modules generated by
$y^\mu_{w,w'}\ot v_{w'.\mu}\ot \xi_{-\nu}$ have larger lowest weight. The corresponding
highest weight with respect to $U_i$ is 
$\nu-w.\mu+(w.\mu,\alpha_i^\vee)\alpha_i=w.\mu-\nu+(\nu,\alpha_i^\vee)\alpha_i$ and hence
\begin{align*}
  E_i^{(\nu,\alpha_i^\vee)+1} y^\mu_{w,w'}\ot (v_{w'.\mu}\ot \xi_{-\nu})=0.
\end{align*}
In view of (\ref{modul-rels}) this proves that $\theta_1$ is well defined. The injectivity
of $\theta_1$ will follow from (\ref{phi1(votxi)}).

To prove (\ref{phi1(votxi)}) note that for any weight vectors $v\in M(w.\mu)_\lambda$ and
$\xi\in M(\nu)^\ast_{\lambda'}$ there exist $E\in \uqlsp_{\lambda'+\nu}$ and 
$F\in \uqlsm_{w.\mu-\lambda}$ such that $v=F v_{w.\mu}$ and $\xi=E \xi_{-\nu}$.
Moreover, in the $\uqls$-module $M(w.\mu)\ot M(\nu)^\ast$ one has
\begin{align*}
  v\ot \xi \in FE(v_{w.\mu}\ot \xi_{-\nu})+\sum_{\zeta < \wght(\xi)}
  M(w.\mu)\ot M(\nu)^\ast_\zeta.
\end{align*}
Hence one obtains by induction on $\wght(\xi)$
\begin{align*}
  \theta_1(v\ot \xi)\in & FE y^\mu_{w,w'}\ot (v_{w'.\mu}\ot \xi_{-\nu})+\sum_{\zeta < \wght(\xi)} 
    V_-\ot M(w.\mu)\ot M(\nu)^\ast_\zeta\\
  = & F y^\mu_{w,w'}\ot (v_{w'.\mu}\ot \xi)+\sum_{\zeta < \wght(\xi)} 
    V_-\ot M(w.\mu)\ot M(\nu)^\ast_\zeta\\
  =& \theta_2(v)\ot\xi + \sum_{\zeta < \wght(\xi)} V_-\ot M(w.\mu)\ot M(\nu)^\ast_\zeta.
\end{align*}
The well definedness and the injectivity of $\thetab_1$ follow from the corresponding 
properties of $\theta_1$ and the relations
\begin{align*}
  (M(w'.\mu)\ot M(\nu)^\ast)_\eta&\cong M(\nu)\ot M(w'.\mu)^\ast\\
  \left(\uqg\ot_{\uqls}(M(w'.\mu)\ot M(\nu)^\ast)\right)_\eta&
    \cong \uqg\ot_{\uqls}(M(\nu)\ot M(w'.\mu)^\ast).
\end{align*}
Formula (\ref{phib1(votxi)}) is proved in the same manner as (\ref{phi1(votxi)}).
$\blacksquare$

\vspace{.5cm}

For any $\mu, \nu\in \psp$ we define
\begin{align*}
  W(\mu,\nu):=\uqg\ot_{\uqls}(M(\mu)\ot M(\nu)^\ast).
\end{align*}
Using the isomorphism 
\begin{align*}
  W(\mu,\nu)\cong V_+V_-\ot M(\mu)\ot M(\nu)^\ast
  \cong V_-V_+\ot M(\mu)\ot M(\nu)^\ast
\end{align*}  
we define two filtration on $W(\mu,\nu)$ as follows
\begin{align}
  \cF_1^k W(\mu,\nu)= \Lin_\C\{u \ot v \ot \xi\,|\, u\in V_-V_+,\,& v\in M(\mu)_\lambda,\label{fil1}\\
     \,&\xi \in M(\nu)^\ast,
     \, \hght (\mu-\lambda)\le k\},\nonumber
\end{align}
\begin{align}
  \cF^k_2 W(\mu,\nu) =\Lin_\C\{u \ot v \ot \xi\,|\, u\in V_+V_-,\,& v\in M(\mu),\label{fil2}\\
  \, &\xi \in M(\nu)^\ast_\lambda,
  \,\hght(\lambda+\nu)\le k\}.\nonumber
\end{align}
\begin{cor}\label{uqls-klein}
  Assume that $w,w'\in W^S$, $w\rightarrow w'$, $\mu\in \wlat^+$, $\nu\in \psp$. 
  The following relation holds in $W(w'.\mu,\nu)$
  \begin{align*}
    \uqg y^\mu_{w,w'}\ot (v_{w'.\mu}\ot \xi_{-\nu}) &\cap \cF_2^k W(w'.\mu,\nu)\\
    &\subseteq \sum_{\hght(\beta)\le k} 
              V_+\uqnm {\uqlsp}_\beta y^\mu_{w,w'}\ot (v_{w'.\mu}\ot \xi_{-\nu}). 
  \end{align*}
  Similarly one has in $W(\nu, w'.\mu)$ the relation
  \begin{align*}
    \uqg x^\mu_{w,w'}\ot (v_\nu\ot \xi_{-w'.\mu}) &\cap \cF_1^k W(w'.\mu,\nu)\\
    &\subseteq \sum_{\hght(\beta)\le k} 
              V_-\uqnp {\uqlsm}_{-\beta} x^\mu_{w,w'}\ot (v_{\nu}\ot \xi_{-w'.\mu}). 
  \end{align*}
\end{cor}
\noindent {\bf Proof:}
Proposition \ref{triangProp} implies the following equalities.
\begin{align*}
  V_+V_-\theta_1(M(w.\mu)\ot M(\nu)^\ast)&\cap \cF_2^k W(w'.\mu,\nu)\\
  \stackrel{(\ref{phi1(votxi)})}{=}& V_+V_-\theta_1\left(\sum_{\hght(\alpha+\nu)\le k} 
                     M(w.\mu)\ot M(\nu)^\ast_\alpha\right) \\
  =& V_+V_-\theta_1\left(\sum_{\hght(\beta)\le k}
         \uqlsm {\uqlsp}_\beta (v_{w.\mu}\ot \xi_{-\nu})\right)\\
  =& \sum_{\hght(\beta)\le k} V_+\uqnm {\uqlsp}_\beta y^\mu_{w,w'}\ot(v_{w.\mu}\ot \xi_{-\nu}).      
\end{align*}
The second relation is verified analogously. 
$\blacksquare$
\begin{cor}\label{notRand}
  Assume that $w,w'\in W^S$, $w\rightarrow w'$, $\nu\in \psp$, and 
  $x\in \uqnp$. Then in $W(w'.0,\nu)$ the relation
  \begin{align}\label{notin-rel1}
    [y_{w,w'},x]\ot(v_{w'.0}\ot \xi_{-\nu})\notin \sum_{\makebox[0cm]{$w''\in W^S,\atop
         w''\rightarrow w'$}}\uqg y_{w'',w'}\ot (v_{w'.0}\ot \xi_{-\nu}) \setminus \{0\}
  \end{align}
  holds. Similarly, for $y\in\uqnm$ the relation
  \begin{align}\label{notin-rel2}
      [x_{w,w'},y]\ot(v_{\nu}\ot \xi_{-w'.0})\notin \sum_{\makebox[0cm]{$w''\in W^S,\atop
           w''\rightarrow w'$}}\uqg x_{w'',w'}\ot (v_{\nu}\ot \xi_{-w'.0}) \setminus \{0\}
  \end{align}
  holds in $W(\nu,w'.0)$
\end{cor}
\noindent {\bf Proof:}
Recall that $y_{w,w'}\in \uqnm_{w.0-w'.0}$. Hence with respect to the decomposition
\begin{align}\label{uqg-decomp}
  \uqg\cong V_+\ot \uqnm\ot \uqlsp\ot U^0
\end{align}
one obtains using (\ref{relations})
\begin{align*}
  [y_{w,w'},x]\in \sum_{\beta>w.0-w'.0} V_+\ot \uqnm_\beta\ot \uqlsp\ot U^0.
\end{align*}
As $U^0$ acts diagonally and $v_{w'.0}$ is a highest weight vector for the action of 
$\uqls$ this implies
\begin{align}\label{est1}
  [y_{w,w'},x]\ot(v_{w'.0}\ot \xi_{-\nu})\in \sum_{\alpha+\beta>w.0} 
    V_+ V_{-,\alpha} \ot M(w'.0)_\beta \ot M(\nu)^\ast.
\end{align}
On the other hand, for any $k\in \N_0$, Corollary \ref{uqls-klein} implies
\begin{align}
  &\sum _{\makebox[1cm]{$w''\in W^S,\atop
         w''\rightarrow w'$}}\uqg y_{w'',w'}\ot (v_{w'.0}\ot \xi_{-\nu}) 
         \cap \cF^k_2W(w'.0,\nu)\nonumber\\
  &\subseteq \sum _{\makebox[1.5cm]{$w''\in W^S,\atop
         w''\rightarrow w'$}}\sum_{\hght(\gamma)\le k} V_+\uqnm{\uqlsp}_\gamma y_{w'',w'}\ot 
      (v_{w'.0}\ot\xi_{-\nu})\nonumber\\      
  &\subseteq \sum _{\makebox[1.5cm]{$w''\in W^S,\atop
         w''\rightarrow w'$}}\sum_{\hght(\gamma+\nu)\le k}
         V_+\uqnm y_{w'',w'}\ot (v_{w'.0}\ot M(\nu)^\ast_\gamma) +\cF_2^{k-1}W(w'.0,\nu)\nonumber\\
  &\subseteq \sum _{\makebox[1.5cm]{$w''\in W^S,\atop
         w''\rightarrow w'$}}\sum_{\alpha+\beta\le w''.0} V_+ V_{-,\alpha}\ot M(w'.0)_\beta
         \ot M(\nu)^\ast +\cF_2^{k-1}W(w'.0,\nu).\label{est2}
\end{align}
Choose now $k\in \N$ such that $[y_{w,w'},x]\ot(v_{w'.0}\ot \xi_{-\nu})\in 
\cF^k_2 W(w'.0,\nu)\setminus \cF^{k-1}_2 W(w'.0,\nu)$ and assume that (\ref{notin-rel1}) 
does not hold. Then (\ref{est1}) and (\ref{est2}) imply that there exists
$w''\in W^S$, $w''\rightarrow w'$ such that $w''.0 >w.0$. This is a contradiction to 
Corollary \ref{uncomparable} 1) and 2). Hence (\ref{notin-rel1}) holds.
Relation (\ref{notin-rel2}) is verified analogously. 
$\blacksquare$

\section{Categorical equivalence}\label{CatEquiv}

{}From now on we will write $\cA=\cqg$ and 
\begin{align}\label{B-def}
  \cB=\{b\in \cA\,|\, b_{(1)}b_{(2)}(k)=\vep(k)b\,
  \mbox{ for all }k\in \uqls\}.
\end{align}
\subsection{Takeuchi's categorical equivalence}\label{Take}
In this subsection Takeuchi's categorical equivalence \cite{a-Tak79}
is recalled in the present special setting. 
Note that $\cB\subset \cA$ is a left coideal subalgebra of the 
Hopf algebra $\cA$. Thus $\Aqr:=\cA/\cB^+\cA$ where $\cB^+=\{b\in \cB\,|\,\vep(b)=0\}$ is a right 
$\cA$-module coalgebra. Moreover, by \cite[Thm.2.2(2)]{a-MullSch} 
$\cA$ is a faithfully flat right $\cB$-module. It was shown in the 
proof of \cite[Thm. 2.2(1),(2)]{a-MullSch} that $\Aqr$ is equal to
the image of $\cA$ under the restriction map 
$\uqg^\circ\rightarrow \uqls^\circ$ of dual Hopf algebras. Therefore the pairing
\begin{align}\label{K-pairing}
  \pair{\cdot}{\cdot}:\uqls\times \Aqr \rightarrow \C
\end{align}
is nondegenerate.
Let $\AqM$ denote the category of finite dimensional left 
$\Aqr$-comodules and let
$\MK$ denote the category of right $\uqls$-modules which are 
isomorphic to  a finite direct sum of modules of the form 
$M(\lambda)^\ast$, $\lambda\in \psp$. The pairing (\ref{K-pairing}) 
induces a functor
\begin{align}\label{Xi-Def}
  \Xi:\AqM\rightarrow \MK
\end{align}
where for $V\in\AqM$ the right $\uqls$-module structure on
$\Xi(V):= V$ is given by $v\ract k=\pair{k}{v_{(-1)}}v_{(0)}$
for all $k\in \uqls$, $v\in V$.
\begin{proposition}\label{Xi-prop}
  The functor $\Xi$ is an equivalence of categories.
\end{proposition}
\noindent {\bf Proof:} By the nondegeneracy of the pairing 
(\ref{K-pairing}) two objects $V,W\in \AqM$ are isomorphic if and only
if the $\uqls$-modules $\Xi(V)$ and $\Xi(W)$ are isomorphic. It remains
to show that all objects of $\MK$ lie in the image of $\Xi$. To this
end, consider the right $\uqls$-module $M(\lambda)^\ast$, 
$\lambda\in P^S_+$. By Lemma \ref{KinU} one can find $\mu\in \wlat_+$
and an embedding of $\uqls$-modules $M(\lambda)\hookrightarrow V(\mu)$.
Then $V(\mu)^\ast$ is a right $\uqg$-module, or equivalently by
definition of $\cA$, a left $\cA$-comodule. Projection onto $\Aqr$
endows $V(\mu)^\ast$ with a left $\Aqr$-comodule structure.
As $V(\mu)$ decomposes into a direct sum of irreducible $\uqls$-modules
the $\uqls$-module $M(\lambda)^\ast$ can be viewed as a direct summand 
of the $\uqls$-module $V(\mu)^\ast$. As the pairing (\ref{K-pairing})
is nondegenerate the $\uqls$-direct summand 
$M(\lambda)^\ast\subset V(\mu)^\ast$ is an $\Aqr$-subcomodule.
By construction, application of $\Xi$ to this $\Aqr$-subcomodule
yields the right $\uqls$-module $M(\lambda)^\ast$.
$\blacksquare$
 
\vspace{.5cm}
Recall that for any coalgebra $C$ the cotensor product of a right
$C$-comodule $P$ and a left $C$-comodule $Q$ is defined by
\begin{align*}
  P\square_C Q :=\left\{\sum_ip_i{\ot} q_i\in P\ot Q\,\bigg|\,\sum_i
   p_{i(0)}{\ot} p_{i(1)}{\ot} q_i=\sum_i p_i{\ot} q_{i(-1)}{\ot}q_{i(0)}
                \right\}.
\end{align*}  
Let $\ABM$ denote the category of left $\cA$-covariant left 
$\cB$-modules. There exist functors
\begin{align}
  \Phi:\ABM&\rightarrow \AqM,& \Phi(\Gamma)&=\Gamma/\B^+\Gamma,\label{Phi-Def}\\
  \Psi:\AqM&\rightarrow \ABM,& \Psi(V)&=\cA\square_{\Aqr} V.\label{Psi-Def}
\end{align}
Here for any $\Gamma\in \ABM$ the left $\Aqr$-comodule structure
on $\Gamma/\B^+\Gamma$ is induced by the left $\cA$-comodule structure 
of $\Gamma$. Moreover, the left $\B$-module and the left
$\cA$-comodule structures of $\cA\square_{\Aqr} V$ are defined on the
first tensor factor.

\begin{theorem}{\em \cite[Theorem 1]{a-Tak79}}\label{catequiv}
  With the notions as above $\Phi$ and $\Psi$ are mutually inverse
  equivalences of categories.
\end{theorem}  

By the above theorem and Proposition \ref{Xi-prop} in order to show
that two $\cA$-covariant $\cB$-modules coincide, it suffices to 
show that the corresponding $\uqls$-modules coincide. This method will
be applied to show that the differential graded algebra
which will be constructed in Section \ref{diffdual} coincides with the
$q$-deformed de Rham complex constructed in \cite{a-heko06}.

A slight refinement of Theorem \ref{catequiv} also takes into account
possible right $\B$-modules structures. 
Let $\ABMB$ and $\AqMB$ denote the categories of left $\cA$-covariant 
$\B$-bimodules and of left $\Aqr$-covariant right $\B$-modules,
respectively. The functors $\Phi$ and $\Psi$ restrict to functors
$\Phi_\B:\ABMB\rightarrow \AqMB$ and $\Psi_\B:\AqMB\rightarrow \ABMB$,
respectively. Here the right $\B$-module structure on
$\Phi_\B(\Gamma)=\Gamma/\B^+\Gamma$ comes from the right $\B$-module
structure of $\Gamma$. The right $\B$-module structure on 
$\Psi_B(V)=\cA\square_{\Aqr} V$ is given by
$(\sum_ip_i\ot q_i)b=\sum_i p_ib_{(-1)}\ot q_i b_{(0)}$.

\begin{cor}\label{catequivcor}
  The functors $\Phi_\B$ and $\Psi_\B$ are mutually inverse equivalences
  of categories.
\end{cor}

\subsection{Locally finite duals of $\uqg$-modules 
induced by $\uqls$-Modules}\label{locallyfinite}

For $\lambda\in \psp$ define
\begin{align*}
  \Om (\lambda):=\{f\in (\WM)^\ast \,|\, \dim (f\uqg)<\infty \}.
\end{align*}
Here the dual vector space $(\WM)^\ast$ of the left 
$\uqg$-module $\WM$ is endowed with a right $\uqg$-module structure 
in the usual way by
$(fu)(v):= f(uv)$ for all $f\in (\WM)^\ast$, $u\in \uqg$, $v\in \WM$.
One has a canonical inclusion
\begin{align*}
  c:\Om (\lambda)\rightarrow \uqg^\ast,\quad f\mapsto 
  c_f:=(u\mapsto f(u\ot v_\lambda)).
\end{align*}
We will freely use the inclusion $c$ to consider $\Om (\lambda)$ as
a subset of $\uqg^*$.

\begin{lemma}
For all $\lambda\in \psp$ one has $\Om (\lambda)\subset \cA$.
In particular one has $\Om (0)=\cB$.
Moreover, $\Om (\lambda)$ is a left $\cA$-covariant $\cB$-bimodule.
\end{lemma}

\noindent 
{\bf Proof: }The dual Hopf algebra $\uqg^\circ$ of $\uqg$ satisfies
\begin{align}\label{Ucirc}
  \uqg^\circ=\{a\in \uqg^\ast\,|\, \dim(a\uqg)<\infty\}.
\end{align}  
Thus by definition $\Om (\lambda)\subset \uqg^\circ$.
Moreover, $\uqg^\circ$ contains $\cA$ as the linear span of the 
matrix coefficients of 
the representations $V(\mu)$, $\mu\in P^+$. Recall that $\uqg$ is 
semisimple and any irreducible finite dimensional representation of $\uqg$ can be 
obtained by tensoring some $V(\mu)$ with a one dimensional 
representation $D_\nu$, $\nu\in\{-1,1\}^r$,
given by  $K_i v=\nu_i v$ for all $v\in D_\nu$. 
As $\lambda\in \psp$ the finite dimensional $\uqg$-module generated
by $c_f$ for $f\in \Om (\lambda)$ decomposes into a direct sum of
irreducible representations isomorphic to $V(\mu)$, $\mu\in \wlat_+$. 
Thus one gets $\Om (\lambda)\subset \cA$.

Note that $\ubar:=W^{M(0)}$ is a left $\uqg$-module coalgebra.
Let $\overline{\phantom{a}}:\uqg\rightarrow \ubar$ denote the 
canonical projection $u\mapsto u\ot v_0$.
Note that $\WM$ is a right and left $\ubar$-comodule, where the
coaction is given by
\begin{align*}
  \kow_L(u\ot v_\lambda)=\overline{u}_{(1)}\ot u_{(2)}\ot v_\lambda
  \in \ubar\ot \WM,\\
  \kow_R(u\ot v_\lambda)=u_{(1)}\ot v_\lambda\ot \overline{u}_{(2)}
  \in \WM\ot\ubar.
\end{align*}
These coactions are compatible with each other and with the 
$\uqg$-module structure of $\WM$. They induce the desired
$\cB$-bimodule structure on $\Om (\lambda)$.
$\blacksquare$

\vspace{.5cm}

\noindent The above lemma implies in particular 
$\Om (\lambda)\in \ABM$. Thus one can apply the functor $\Phi$
from the previous subsection. The following proposition states
that up to dualization $\Om $ is the inverse of $\Xi\circ\Phi$.

\begin{proposition}\label{TakLokFin}
  For all $\lambda\in \psp$ one has 
  \begin{align*}
   \Xi(\Phi(\Om (\lambda)))=M(\lambda)^\ast.
  \end{align*} 
\end{proposition}
\noindent{\bf Proof:}
Note first that by definition of the left $\cB$-module structure
of $\Om (\lambda)$ one has $(\cB^+\Om (\lambda))(1\ot M(\lambda))=0$.
Thus there exists a well defined pairing
\begin{align}\label{p-lambda}
  \pair{\cdot}{\cdot}_\lambda:\,
  \Om (\lambda)/(\cB^+\Om (\lambda))\times M(\lambda)\longrightarrow \C.
\end{align}
The pairing $\pair{\cdot}{\cdot}_\lambda$ induces a map of right
$\uqls$-modules 
\begin{align*}
  \varphi:\Xi(\Phi(\Om (\lambda)))\rightarrow M(\lambda)^\ast.
\end{align*}  
As $\Om (\lambda)\subset \cA$ the induced map of quotients
\begin{align*}
  i:\Om (\lambda)/(\cB^+\Om (\lambda))\rightarrow \cA/\cB^+\cA
\end{align*}
is also injective by Theorem \ref{catequiv}. Moreover, let $\pi$ denote
the surjection
\begin{align*}
  \pi: \uqls\rightarrow M(\lambda),\quad k\mapsto kv_\lambda.
\end{align*}
Then the pairings (\ref{K-pairing}) and (\ref{p-lambda}) satisfy
\begin{align*}
  \pair{\overline{f}}{\pi(k)}_\lambda=\pair{k}{i(\overline{f})}
\end{align*}
for all $\overline{f}\in\Om (\lambda)/(\cB^+\Om (\lambda))$, 
$k\in \uqls$.
As the pairing (\ref{K-pairing}) is nondegenerate and $i$ is injective
this implies that $\varphi$ is injective. By Theorem \ref{catequiv}
as $M(\lambda)$ is irreducible it remains to show that 
$\Om (\lambda)\neq 0$. To this end apply Lemma~\ref{KinU} to pick 
$\mu\in \wlat^+$ such that there exists an embedding 
$M(\lambda)\hookrightarrow V(\mu)$ of $\uqls$-modules. Let $v$ denote
the image of $v_\lambda$ under this embedding. Pick $g\in V(\mu)^\ast$
such that $g(v)\neq 0$ and let $c_{g,v}\in\cqg$ denote the 
corresponding matrix coefficient. Then there is an element
$f\in \Om (\lambda)\setminus \{0\}$ defined by 
$f(u\ot v_\lambda)=c_{g,v}(u)$. $\blacksquare$

\section{$q$-Differential forms as locally finite duals}\label{diffdual}
{From} now on we restrict to the case of irreducible flag manifolds
$G/P_S$. Thus $S=\pi\setminus \{\alpha_s\}$ where $\alpha_s$ occurs
in each positive root of $\gfrak$ with multiplicity at most one.
Let again $\cB\subset \cqg$ be the left coideal subalgebra defined by
(\ref{B-def}).

\subsection{$q$-Differential forms for irreducible flag manifolds}\label{q-Diffforms}
The aim of this section is to recall the structure
of the canonical differential graded algebra over $\B$ constructed 
and investigated in \cite{a-heko06}, \cite{a-heko04}.

To this end recall that a first order differential calculus (FODC)
over $\B$ is a $\B$-bimodule $\Gamma$ together with a $\C$-linear
map
\begin{equation*}
  \dif:\B\rightarrow\Gamma
\end{equation*}
such that $\Gamma=\Lin_\C\{a\,\dif b\,c\,|\,a,b,c\in\B\}$ and $\dif$
satisfies the Leibniz rule
\begin{align*}
  \dif(ab)&=a\,\dif b + \dif a\,b.
\end{align*}    
If $\Gamma$ possesses the structure of a left $\cA$-comodule
\begin{equation*}
  \kow_\Gamma:\Gamma\rightarrow\cA\ot \Gamma
\end{equation*}
such that
\begin{equation*}
\kow_\Gamma(a\dif b\,c)=(\kow_\B a)((\id\otimes\dif)\kow_\B b)
                          (\kow_\B c)
\end{equation*}
then $\Gamma$ is called (left) covariant.
A covariant FODC $\Gamma\neq \{0\}$ over $\B$ is called
irreducible if it does not possess any nontrivial quotient
(by a left covariant $\B$-bimodule). The dimension of a covariant
FODC $\Gamma\neq \{0\}$ over $\B$ is defined by 
$\dim \Gamma=\dim_\C \Gamma/\B^+\Gamma$.
Any finite dimensional covariant FODC over $\B$ is uniquely determined by its so called 
quantum tangent space
\begin{align*}
  T_\Gamma=\{f\in \Gamma^\ast\,|\, f|_{\B^+\Gamma}=0\},
\end{align*}
(see \cite[Lemma 6]{a-HK-QHS}, \cite[Remark 2.4]{a-heko06}).
The quantum tangent space can be considered as a subset of the dual
coalgebra $\B^\circ$ of $\B$ via the map $f\mapsto (b\mapsto f(\dif b))$.
It is one of the main results of \cite{a-heko04} that there exist 
precisely two finite dimensional irreducible covariant FODC $(\Gamma_\del,\del)$ and
$(\Gamma_\delb,\delb)$ over $\B$. The quantum tangent
spaces of the FODC $\Gamma_\del$ and $\Gamma_\delb$ \cite[Propositions 3.3, 3.4]{a-heko06} 
are given by
\begin{align}\label{tangentspaces}
  T_\del= (\ad\,\uqls) F_s,\qquad T_\delb= (\ad\, \uqls) E_s,
\end{align}
respectively, considered as subspaces of $\B^\circ$ via the pairing (\ref{pair}). Moreover, 
the FODC $\Gamma_\del$ and $\Gamma_\delb$ satisfy
\begin{align}\label{B+commute}
  \B^+\Gamma_\del=\Gamma_\del \B^+,\qquad
  \B^+\Gamma_\delb=\Gamma_\delb \B^+.
\end{align}
The direct sum $\Gamma_\dif=\Gamma_\del\oplus\Gamma_\delb$ with 
the map $\dif=\del\oplus\delb$ is a covariant FODC which is
a $q$-analog of the K\"ahler differentials over the affine 
algebraic variety $G/L_S$.

Consult \cite[Section 2.3.2]{a-heko06} for the definition of the universal
differential calculus of a FODC $(\Gamma,\dif)$. Let $\Gduw$,
$\Gdbuw$, and $\Gdifuw$ denote the universal differential calculi
of the FODC $(\Gamma_\del,\del)$, $(\Gamma_\delb,\delb)$, and
$(\Gamma_\dif,\dif)$, respectively. The following theorem is
contained in \cite[Propositions 3.6, 3.7, 3.11]{a-heko06}.
\begin{theorem}\label{summaryTheo}
  \begin{itemize}
    \item[(i)] The multiplicity of weight spaces of the left
    $\uqls$-module 
    $(\Xi\circ\Phi(\Gduw[k]))^\ast=(\Gduw[k]/\B^+\Gduw[k])^\ast$
    coincides with the multiplicity of weight spaces of the left
    $U(\lfrak_S)$-module $\Lambda^k(\gfrak/\pfrak_S)$. In particular
    \begin{align*}
      \dim_\C(\Gduw[k]/\B^+\Gduw[k])={\dim_\C(\gfrak/\pfrak_S)\choose k}.
    \end{align*}
    \item[(ii)] The multiplicity of weight spaces of the left
    $\uqls$-module 
    $(\Xi\circ\Phi(\Gdbuw[k]))^\ast=(\Gdbuw[k]/\B^+\Gdbuw[k])^\ast$
    coincides with the multiplicity of weight spaces of the left
    $U(\lfrak_S)$-module $\Lambda^k(\gfrak/\pfrak_S)^\ast$. In particular
    \begin{align*}
      \dim_\C(\Gdbuw[k]/\B^+\Gdbuw[k])={\dim_\C(\gfrak/\pfrak_S)\choose k}.
    \end{align*}
    \item[(iii)] For all $k\in \N_0$ the canonical map
      \begin{align}\label{isom}
         \bigoplus_{i+j=k}\Gduw[i]/\B^+\Gduw[i]\ot\Gdbuw[j]/\B^+\Gdbuw[j]
                       \rightarrow \Gdifuw[k]/\B^+\Gdifuw[k]
      \end{align}
      is an isomorphism. In particular
      \begin{align*}
        \dim \Gdifuw[k]={2\dim_\C(\gfrak/\pfrak_S)\choose k}.
      \end{align*}  
   \end{itemize}
\end{theorem}
 The above theorem implies in particular, that 
 $\Phi(\Gdifuw[2 \dim(\gfrak/\pfrak_S)])$ is the trivial
 one-dimensional left $\Aqr$-comodule. Moreover, by (\ref{B+commute})
 the right $\B$-action on $\Phi(\Gdifuw[2 \dim(\gfrak/\pfrak_S)])$
 is trivial, i.e.~$\gamma b=\vep(b)\gamma$ for all $b\in \B$, 
 $\gamma\in \Phi(\Gdifuw[2 \dim(\gfrak/\pfrak_S)])$.
 Hence by the categorical equivalence in Corollary \ref{catequivcor} the 
 covariant $\B$-bimodules $\Gdifuw[2 \dim(\gfrak/\pfrak_S)]$ and 
 $\B$ are isomorphic. This observation implies the following corollary.
 \begin{cor}
 $\Gdifuw[2 \dim(\gfrak/\pfrak_S)]$
 is a free left and right $\B$-module generated by one left 
 $\cqg$-coinvariant element 
 $\omvol\in \Gdifuw[2 \dim(\gfrak/\pfrak_S)]$ satisfying
 $\omvol \,b=b\,\omvol$ for all $b\in \B$.
 \end{cor}

\subsection{The differential calculus $\Gduw$}\label{delcalc}
One is now in a position to construct the differential graded algebras $\Gduw$, $\Gdbuw$,
and $\Gdifuw$ from \cite{a-heko06} as locally finite duals of BGG-like sequences
of $\uqg$-modules induced by $\uqls$-modules. We begin with $\Gduw$.
Consider the BGG-resolution 
 \begin{align}\label{qBGGV(0)}
   C^S_{\ast,0}:\qquad 0\longrightarrow C_{\dim \gfrak/\pfrak_S,0}^S
    \stackrel{\varphi^S_{\dim(\gfrak/\pfrak_S)}}{\longrightarrow} \cdots 
    \stackrel{\varphi^S_1}{\longrightarrow} C_{0,0}^S
    \stackrel{\vep_\mu}{\longrightarrow} V(0)\longrightarrow 0,
  \end{align}
of the trivial $\uqg$-module $V(0)$, the corresponding sequence (\ref{CundSeq}) obtained
by applying the functor $\fun$
  \begin{align*}
   \Cu^S_{\ast,0}:\, 0\longrightarrow 
   \Cu_{\dim \gfrak/\pfrak_S,0}^S
    \stackrel{\vphiu^S_{\dim(\gfrak/\pfrak_S)}}{\longrightarrow} 
    \cdots 
    \stackrel{\vphiu^S_1}{\longrightarrow} \Cu_{0,0}^S
    \stackrel{\underline{\vep}_\mu}{\longrightarrow} \uqg\ot_{\uqpsop}V(0)
    \longrightarrow 0,
  \end{align*}
and its locally finite dual
\begin{align*}
   \Om ^{\ast,0}:\qquad0 \longleftarrow \Om ^{\dim \gfrak/\pfrak_S,0}
   \stackrel{\del_{\dim(\gfrak/\pfrak_S)}}{\longleftarrow} \cdots 
    \stackrel{\del_1}{\longleftarrow} \Om ^{0,0}\cong \cB
    \stackrel{}{\longleftarrow} \C\longleftarrow 0,
  \end{align*}
where 
\begin{align*} 
  \Om ^{n,0}\cong \bigoplus_{w\in W^S, l(w)=n}\Om (w.0).
\end{align*}  
Recall from Section \ref{BGGsection} that the differentials of the complexes
$C^S_{\ast,0}$ and $\Cu^S_{\ast,0}$  are given in terms of 
elements $y_{w,w'}\in \uqnm$ where $w,w'\in W^S$ and $w\rightarrow w'$.
In the case of irreducible flag manifolds the simple reflection
$s_s$ corresponding to $\alpha_s$ is the only element in $W^S$
of length one. Note that $s_s.0=-\alpha_s$. Thus the differential
$\vphiu_1:V^{M(-\alpha_s)}\rightarrow V^{M(0)}$
is determined by $y_{s_s,e}=F_s$ up to multiplication by a nonzero
factor. The corresponding differential
\begin{align*}
  \del_1:\cB\cong\Om ^{0,0}\rightarrow \Om (-\alpha_s)\cong\Om ^{1,0}
\end{align*}
satisfies the Leibniz rule.
Indeed, for all $a,b\in \cB$, $u\in \uqg$ one has
\begin{align*}
  \del_1(ab)(u\ot v_{-\alpha_s})&=(ab)(u F_s\ot v_0)\\
  &=a(u_{(1)}F_s\ot v_0) b(u_{(2)}\ot v_0)+ 
   a(u_{(1)}K^{-1}_s\ot v_0) b(u_{(2)}F_s\ot v_0)\\
  &=(\del_1(a)b +a\del_1(b))(u\ot v_{-\alpha_s}).
\end{align*}
\begin{lemma}\label{FODClemma}
  $(\del_1{:}\cB\rightarrow \Om ^{1,0})$ is a covariant FODC isomorphic to 
  $(\del{:}\cB\rightarrow \Gamma_\del)$.
\end{lemma}
\noindent 
{\bf Proof:} Recall from (\ref{tangentspaces}) that $T_\del$ is an irreducible $\uqls$-module
of highest weight $-\alpha_s$. Taking duals one obtains that $M(-\alpha_s)^\ast$ is 
isomorphic to $\Gamma_\del/\B^+\Gamma_\del$ as a right $\uqls$-module. Proposition
\ref{TakLokFin} and the categorical equivalence now imply that 
$\Om^{1,0}\cong\Om (-\alpha_s)$ is an  $\cA$-covariant
$\cB$-bimodule isomorphic to $\Gamma_\del$. 

As $M(-\alpha_s)$ is an irreducible $\uqls$-module it remains to check 
that $\del_1\neq0$. This is a special case of the following Lemma
which is proved independently of the above claim. 
$\blacksquare$

\begin{lemma}\label{surlemma}
  For all $n\in \N_0$ the map
  \begin{align*}
    \psi_n:\cB\ot \Om ^{n,0}\rightarrow \Om ^{n+1,0},\qquad 
        b\ot \omega\mapsto b\,\del_{n+1}\omega
  \end{align*}
  is surjective.
\end{lemma}
\noindent {\bf Proof:} It suffices to show that for any 
$w,w'\in W^S$ such that $w\rightarrow w'$ one has
\begin{align}\label{yww'not0}
  y_{w,w'}\notin \uqg \uqls^+
\end{align}
where $\uqls^+=\ker \vep\cap \uqls$ denotes the augmentation ideal of $\uqls$.
Indeed, choose $f\in \Om (w'.0)$ such that $f(1\ot v_{w'.0})\neq 0$.
Then for all $b\in \cB$ one has
\begin{align*}
  (\del_{n+1}(bf))(1\ot v_{w.0})=
  b(y_{w,w'(1)})f(y_{w,w'(2)}\ot v_{w'.0}).
\end{align*}
Since $\cB$ separates $\ubar=\uqg/\uqg \uqls^+$ 
\cite[Prop.~6.1]{a-heko04} 
and $y_{w,w'}\notin \uqg \uqls^+$ by assumption (\ref{yww'not0}) and
$y_{w,w'}\in \uqnm$, one can choose $b\in \cB$ such that
\begin{align*}
  b(y_{w,w'(1)}) y_{w,w'(2)}=1.
\end{align*}
By assumption on $f$ this implies 
\begin{align}\label{delbfneq0}
  \del_{n+1}(bf)|_{W^{M(w.0)}}\neq 0.
\end{align}
By the categorical equivalence and Proposition \ref{TakLokFin} the covariant
$\cB$-module $\Om ^{n+1,0}$ contains any irreducible covariant
$\cB$-submodule with multiplicity at most one.
Since $\im \psi_n$ is a covariant left $\B$-module relation (\ref{delbfneq0})
implies $\Om (w.0)\subset \im \psi_n$.

It remains to verify (\ref{yww'not0}). Assume on the contrary that
$y_{w,w'}\in V_-\uqlsm^+$ where $\uqlsm^+=\ker \vep\cap\uqlsm$.
Then $y_{w,w'}\ot v_{w'.0}\in V_-\ot \uqlsm^+v_{w'.0}\subset V^{M(w'.0)}$ is nonzero 
because by Remark \ref{standnonzero} the standard map does not vanish.
Since $(\ad \uqls)V^-=V^-$ and $M(w'.0)$ is an irreducible
$\uqls$-module there exists $E_i$, where $i\neq s$, such that 
$E_i y_{w,w'}\ot v_{w'.0}\neq 0$.
This is a contradiction to $h_{w,w'}(E_i\ot v_{w.0})=0$. $\blacksquare$

\vspace{.5cm}

By Lemma \ref{FODClemma} one has $\Om^{1,0}=\Lin_\C\{a\,\del_1 b\,|\, a,b\in \B\}$. 
Define a map
\begin{align}\label{barwedgedef}
  \barwedge:\Om ^{1,0}\ot_\cB \Om ^{n,0}&\rightarrow \Om ^{n+1,0}\\
  (a\,\del_1 b)\ot \omega &\mapsto a\,\del_1 b\barwedge \omega
    :=a\,\del_{n+1}(b\omega) - ab\,\del_{n+1}\omega.\nonumber
\end{align}
\begin{lemma}\label{barwedgeWellDef}
  The map $\barwedge$ is well defined.
\end{lemma}
\noindent
{\bf Proof:} Recall from the last statement of Corollary 
\ref{uncomparable} that $y_{w,w'}\in \uqnm_{-\beta}$
where $\omega_s(\beta)=1$. Therefore 
\begin{align}\label{ob-ot-id}
  (\ob\ot\id)\kow y_{w,w'} - 1\ot y_{w,w'}\in\overline{\uqls F_s}\ot \uqls \subset\ubar\ot \uqg
\end{align}
where as before $\ubar=\uqg/\uqg\uqls^+$.
To prove that $\barwedge$ is well defined consider 
$a_i,b_i\in \cB$ such that $\sum_i a_i\,\del_1b_i=0$ or equivalently
\begin{align}\label{assumption}
  \sum_i a_i(u_{(1)}) b_i(u_{(2)}F_s)=0\qquad \forall u\in \uqg.
\end{align}
Since $a_i(ux)=a_i(u)\vep(x)$ for all $u\in \uqg$, $x\in \uqls$ formula (\ref{assumption})
is equivalent to
\begin{align}\label{mit-x}
  \sum _i a_i(u_{(1)}) b_i(u_{(2)} x F_s)=0\qquad \forall u\in \uqg, x\in \uqls.
\end{align}
Then for all $\omega\in \Om ^{n,0}$ and all $u\in \uqg$ one has
\begin{align*}
  \sum_i(a_i\del_{n+1}&(b_i\omega)-a_ib_i\del_{n+1}\omega)
  (u\ot v_{w.0})
  =\sum_i a_i(u_{(1)})
    \Bigg[\\
    &(b_i\omega)\bigg(\sum_{w'}u_{(2)}y_{w,w'}\ot v_{w'.0}\bigg)
   -b_i(u_{(2)})\omega\bigg(\sum_{w'}u_{(3)}y_{w,w'}\ot 
                   v_{w'.0}\bigg)\Bigg]\\
  =&\sum_i a_i(u_{(1)})\sum_{w'}b_i(u_{(2)}y_{w,w'(1)}^+)
     \omega(u_{(3)}y_{w,w'(2)}\ot v_{w'.0})=0
\end{align*}
where $y^+=y-\vep(y)$ and the last equation follows from (\ref{ob-ot-id}) and (\ref{mit-x}).
Thus $\barwedge:\Om ^{1,0}\times \Om ^{n,0}\rightarrow \Om ^{n+1,0}$ 
is well defined. Moreover, by definition for 
$a,b\in \cB$ and $\omega\in\Om ^{n,0}$ one has
\begin{align*}
  ((\del_1 b)a)\barwedge \omega=(\del_1(ba)-b\del_1a)\barwedge\omega
  =\del_{n+1}(ba\omega) - b\del_{n+1}(a\omega)
  = \del_1 b\barwedge a\omega
\end{align*}
and thus $\barwedge$ is defined on $\Om ^{1,0}\ot_\cB \Om ^{n,0}$. 
$\blacksquare$
\begin{lemma}\label{barwedgeSur}
  The map $\barwedge:\Om ^{1,0}\ot_\cB \Om ^{n,0}\rightarrow 
  \Om ^{n+1,0}$ is surjective.
\end{lemma}
\noindent{\bf Proof:} 
One shows by induction that for all $a_1,\dots,a_k\in \B$ the relation
\begin{align}\label{boundary}
  \del_1 a_1\barwedge(\del_1 a_2\barwedge(\dots
         \barwedge \del a_k)\dots))=
  \del_k( a_1 \del_1 a_2\barwedge(\dots
         \barwedge \del a_k)\dots))     
\end{align}
holds. The claim of the lemma holds for $n=0$. Using Lemma \ref{surlemma} and 
(\ref{boundary}) one shows by induction on $n$ that
\begin{align*}
  \Om ^{n,0}=\Lin_\C\{a_0\del_1 a_1\barwedge(\del_1 a_2\barwedge(\dots
         (\barwedge \del_1 a_n)\dots))\,|\, a_0,a_1,\dots,a_n\in\cB\}.\qquad \blacksquare
\end{align*}

\vspace{.5cm}

As in \cite{a-heko06} let $\Gduw$ denote the universal differential
calculus with FODC $(\del:\cB\rightarrow \Gamma_\del)$.
Define a map 
\begin{align*}
  \iota^n:\Gduw[n]&\rightarrow \Om ^{n,0}\\
   a_0\del a_1\wedge\del a_2
  \wedge \dots \wedge \del a_n &\mapsto a_0 \del_1 a_1\barwedge
  (\del_1 a_2\barwedge( \dots \barwedge \del_1 a_n))
\end{align*}  
by repeated use of the map $\barwedge$.
\begin{lemma}\label{ikWellDef}
  The map $\iota^n$ is well defined.
\end{lemma}
\noindent {\bf Proof:} 
By definition of $\Gduw$ it suffices to show that
for all $a_i,b_i\in \cB$ such that $\sum_i a_i\del b_i=0$ and for all
$\omega\in \Om ^{k,0}$, where $0\le k\le n-2$, one has 
$\sum_i\del_1 a_i\barwedge (\del_1 b_i\barwedge \omega)=0$.
This can be seen as follows.
\begin{align*}
  \sum_i &\del_1 a_i\barwedge (\del_1 b_i\barwedge \omega)
  =\sum_i \del_{k+2}(a_i\del_1 b_i\barwedge \omega)
    -\sum_i a_i \del_{k+2}(\del_{1} b_i\barwedge \omega)\\
  &=-\sum_i a_i\del_{k+2}(\del_{k+1}(b_i\omega)-b_i\del_{k+1}\omega)
  =\sum_i a_i\del_1 b_i\barwedge \del_{k+1 }\omega=0.\,\blacksquare
\end{align*}

Note that by construction $\iota^\ast$ is a morphism of complexes.
Moreover, by Theorem \ref{summaryTheo}(i) and Proposition \ref{C=D}
the dimensions of the covariant left $\cB$-bimodules $\Gduw[n]$ and 
$\Om ^{n,0}$ coincide. As $\iota^n$ is a surjective map of covariant 
left $\cB$-modules by Lemma \ref{barwedgeSur} the categorical 
equivalence implies that $\iota^n$ is an isomorphism.
\begin{proposition}\label{DelCalcIso}
  The map $\iota^\ast:\Gduw[\ast] \rightarrow \Om ^{\ast,0}$ is an 
  isomorphism of complexes of covariant $\cB$-bimodules.
\end{proposition}

\subsection{The differential calculus $\Gdbuw$} \label{delbcalc}
Recall from Section \ref{q-Diffforms} that there exists a second irreducible covariant
FODC $(\Gamma_\delb,\delb)$ over $\cB$. It follows from (\ref{tangentspaces}) that
$T_{\delb}\cong M(-\alpha_s)^\ast$. As in the case $\Gduw$ the universal differential
calculus $\Gdbuw$ can be obtained as the 
locally finite dual of a sequence of $\uqg$-modules induced by 
$\uqls$-modules.
This can be seen using the involutive algebra automorphism
coalgebra antiautomorphism $\invol:\uqg\rightarrow \uqg$ defined in Section \ref{uqgcqg}.
The exact sequences $C_{\ast,0}^S$ and $\Cu_{\ast,0}^S$ from the 
previous subsection can be endowed with a new $\uqg$-module structure
via $\invol$. Using the isomorphism 
$V_\invol^{M(\lambda)}\cong V^{M(\lambda)^\ast}$ and
$W_\invol^{M(\lambda)}\cong W^{M(\lambda)^\ast}$ one obtains complexes
\begin{align*}
   C^S_{0,\ast}:\qquad 
   0\longrightarrow C_{0,\dim \gfrak/\pfrak_S}^S
    \stackrel{\ophi^S_{\dim(\gfrak/\pfrak_S)}}{\longrightarrow} \cdots 
    \stackrel{\ophi^S_1}{\longrightarrow} C_{0,0}^S
    \stackrel{\oep_\mu}{\longrightarrow} V(0)\longrightarrow 0,
  \end{align*}
and
  \begin{align*}
   \Cu^S_{0,\ast}:\, 0\longrightarrow 
   \Cu_{0,\dim \gfrak/\pfrak_S}^S
    \stackrel{\underline{\ophi}^S_{\dim(\gfrak/\pfrak_S)}}{\longrightarrow} 
    \cdots 
    \stackrel{\underline{\ophi}^S_1}{\longrightarrow} \Cu_{0,0}^S
    \stackrel{\underline{\oep}_\mu}{\longrightarrow} \uqg\ot_{\uqps} V(0)
    \longrightarrow 0,
  \end{align*}
where 
\begin{align*}
  C^S_{0,n}=\bigoplus_{w\in W^S,\, l(w)=n}V^{M(w.0)^\ast},\qquad
  \Cu^S_{0,n}=\bigoplus_{w\in W^S,\, l(w)=n}W^{M(w.0)^\ast}. 
\end{align*}
If $w,w'\in W^S$ and $w\rightarrow w'$ then the component of the 
differential which maps $V^{M(w.0)^\ast}\rightarrow V^{M(w'.0)^\ast}$
is given by
\begin{align*}
  u\ot \xi_{-w.0}\mapsto u x_{w,w'}\ot \xi_{-w'.0}
\end{align*}
where $x_{w,w'}=\eta(y_{w,w'})$.
Taking locally finite duals one obtains a complex
\begin{align*}
   \Om ^{0,\ast}:\qquad\Om ^{0,\dim \gfrak/\pfrak_S}
   \stackrel{\delb_{\dim(\gfrak/\pfrak_S)}}{\longleftarrow} \cdots 
    \stackrel{\delb_1}{\longleftarrow} \Om ^{0,0}\cong \cB
    \stackrel{}{\longleftarrow} \C\longleftarrow 0.
  \end{align*}
To show that $\Om ^{0,\ast}$ is isomorphic as a complex of covariant
$\cB$-bimodules to the complex $\Gdbuw$ constructed in 
\cite[3.3.2]{a-heko06} the arguments of the last subsection can be 
repeated.

\subsection{The differential calculus $\Gdifuw$}\label{difcalc}
Combining the constructions from the previous two subsections
the $q$-analog of the de Rham-complex over $G/L_S$ can also be 
realized as a locally finite dual of a sequence of $\uqg$-modules
induced by $\uqls$-modules.
To this end define
\begin{align}\label{Ctilde_nm}
  \Cu_{n,m}^S:=\bigoplus_{\makebox[0cm]{$w_1,w_2\in W^S\atop
              l(w_1)=n,\, l(w_2)=m$}}W(w_1.0,w_2.0).
\end{align}              
Recall that for each $w_1,w_2\in W^S$ the $\uqg$-module
$W(w_1.0,w_2.0)$ is a cyclic module generated by 
$1\ot (v_{w_1.0}\ot \xi_{-w_2.0})$.
If $w_1,w_1',w_2\in W^S$ and $w_1\rightarrow w_1'$ define a map
\begin{align*}
  \hu_{w_1,w_1';w_2}:W(w_1.0,w_2.0)&\rightarrow W(w_1'.0,w_2.0),\\
  u\ot (v_{w_1.0}\ot \xi_{-w_2.0})&\mapsto 
  uy_{w_1,w_1'}\ot (v_{w'_1.0}\ot \xi_{-w_2.0}).
\end{align*}
Similarly, if $w_1,w_2,w'_2\in W^S$ and $w_2\rightarrow w_2'$ define 
a map
\begin{align*}
  \hu_{w_1;w_2,w'_2}:W(w_1.0,w_2.0)&\rightarrow W(w_1.0,w'_2.0),\\
  u\ot (v_{w_1.0}\ot \xi_{-w_2.0})&\mapsto 
  ux_{w_2,w_2'}\ot (v_{w_1.0}\ot \xi_{-w'_2.0}).
\end{align*}
Note that the symbol $\fun$ in the above definitions of $\Cu_{n,m}^S$, $\hu_{w_1,w_1';w_2}$,
and $\hu_{w_1,w_1';w_2}$
is only a formal reminiscence of the functor $\fun$ from Section \ref{functor}.
No functorial properties will be used in the present section. 
\begin{lemma}
  For all $w_0,w_1,w_2\in W^S$ such that $w_1\rightarrow w_2$
  the maps $\hu_{w_1,w_2;w_0}$ and $\hu_{w_0;w_1,w_2}$ are well defined.
\end{lemma}
\noindent{\bf Proof:} Note that 
\begin{align*} 
  \uqg\ann_{\uqls}(v_{w_1.0}\ot \xi_{-w_0.0})=
  \uqg\ann_{\uqls}(y_{w_1,w_2}\ot (v_{w_2.0}\ot \xi_{-w_0.0}))
\end{align*} 
because $M(w_1.0)\ot M(w_0.0)^\ast$ and 
$\uqls y_{w_1,w_2}\ot (v_{w_2.0}\ot \xi_{-w_0.0})\subset W(w_2.0,w_0.0)$
are isomorphic as $\uqls$-modules by Proposition \ref{triangProp}.
This proves that $\hu_{w_1,w_2;w_0}$ is well defined.
The second statement follows analogously.
$\blacksquare$.

\vspace{.5cm}

For each $w_2\in W^S$ one obtains a sequence
\begin{align}\label{astw2sequence}
   \Cu^S_{\ast,w_2}:\qquad 
     \dots
    \stackrel{\hu_{n+1,w_2}}{\longrightarrow} \Cu^S_{n,w_2}
    \stackrel{\hu_{n,w_2}}{\longrightarrow} 
    \Cu^S_{n-1,w_2}
    \stackrel{\hu_{n-1,w_2}}{\longrightarrow} \dots
    \stackrel{\hu_{1,w_2}}{\longrightarrow}\Cu^S_{0,w_2}
\end{align}
where 
\begin{align*}
  \Cu^S_{n,w_2}=\bigoplus_{w_1\in W^S,\, l(w_1)=n}W(w_1.0,w_2.0),
  \qquad\hu_{n,w_2}=
         \sum_{\makebox[0cm]{$w_1,w_1'\in W^S,\, l(w_1)=n\atop
         w_1\rightarrow w_1'$}}\hu_{w_1,w_1';w_2}.
\end{align*}              
This sequence satisfies $\hu_{n,w_2}\hu_{n+1,w_2}=0$ for all $n\in \N$.
Indeed, the exactness of the sequence (\ref{CundSeq}) implies that for 
any $w_1,w_1''\in W^S$ such that $l(w_1)=n+1$ and $l(w''_1)=n-1$ one 
has
\begin{align*}
  \sum_{\makebox[0cm]{$w_1'\in W^S, \atop
         w_1\rightarrow w_1'\rightarrow w''_1$}}
         y_{w_1,w'_1}y_{w_1',w''_1}\in \uqnm \ann_{\uqlsm}(v_{w''_1.0}).
\end{align*}

Similarly, for each $w_1\in W^S$ one has a complex
\begin{align}\label{w1astsequence}
   \Cu^S_{w_1,\ast}:\qquad 
     \dots
    \stackrel{\hu_{w_1,n+1}}{\longrightarrow} \Cu^S_{w_1,n}
    \stackrel{\hu_{w_1,n}}{\longrightarrow} 
    \Cu^S_{w_1,n-1}
    \stackrel{\hu_{w_1,n-1}}{\longrightarrow} \dots
    \stackrel{\hu_{w_1,1}}{\longrightarrow}\Cu^S_{w_1,0}
\end{align}
where 
\begin{align*}
  \Cu^S_{w_1,n}=\bigoplus_{w_2\in W^S,\, l(w_2)=n}W(w_1.0,w_2.0),
  \qquad  \hu_{w_1,n}=
         \sum_{\makebox[0cm]{$w_2,w_2'\in W^S,\, l(w_2)=n\atop
         w_2\rightarrow w_2'$}}\hu_{w_1;w_2,w_2'}.
\end{align*}    
To prove exactness of the sequences (\ref{astw2sequence}) and (\ref{w1astsequence})
we extend the filtrations defined in Section \ref{Homest} to the vector spaces 
$\Cu^S_{n,w_2}$ and $\Cu^S_{w_1,n}$.
Define 
\begin{align*}
  \cF^k_2 \Cu^S_{n,w_2} = \bigoplus_{\makebox[1cm]{$w_1\in W^S,\atop l(w_1)=n$}} 
                                        \cF^k_2 W(w_1.0,w_2.0),\qquad
  \cF^k_1 \Cu^S_{w_1,n} = \bigoplus_{\makebox[1cm]{$w_2\in W^S,\atop l(w_2)=n$}} 
                                        \cF^k_1 W(w_1.0,w_2.0).
\end{align*}
\begin{lemma}\label{filtered}
  For any $w_1,w_2\in W^S$ the complexes $\Cu^S_{w_1,\ast}$ and $\Cu^S_{\ast,w_2}$ are filtered
  by the filtrations $\cF_1$ and $\cF_2$, respectively.
\end{lemma}

\noindent{\bf Proof:} 
Consider $w_1,w_1', w_2'\in W^S$ such that $w_1\rightarrow w_1'$. We show that
\begin{align*}
  \hu_{w_1,w_1';w_2}(\cF^k_2 W(w_1.0,w_2.0))\subset \cF^k_2 W(w_1'.0,w_2.0).
\end{align*}  
Define $\cF^k M(w_2.0)^\ast=\bigoplus_{\hght(\mu+w_2.0)\le k} 
  M(w_2.0)^\ast_\mu$. Then
\begin{align*}
  \cF^k_2 W(w_1.0,w_2.0)&\subset V_+\uqnm \ot v_{w_1.0} \ot \cF^k M(w_2.0)^\ast\\
           &=\sum_{\hght(\beta)\le k} V_+\uqnm \uqlsp_\beta\ot (v_{w_1.0} \ot \xi_{-w_2.0}).
\end{align*}
Therefore
\begin{align*}
  \hu_{w_1,w_1';w_2}(\cF^k_2 W(w_1.0,w_2.0))&\subset \sum_{\hght(\beta)\le k}
                         V_+\uqnm\uqlsp_\beta y_{w_1,w_1'}
                          \ot (v_{w_1'.0}\ot \xi_{-w_2.0})\\
              & \subset V_+\uqnm \ot (v_{w_1'.0}\ot \cF^k M(w_2.0)^\ast)\\           
              & \subset V_+V_- \ot M(w_1'.0)\ot \cF^k M(w_2.0)^\ast.
\end{align*}
The statement about $\Cu^S_{w_1,\ast}$ and $\cF_1$ is verified analogously. $\blacksquare$

\medskip

Let $\gr_{\cF_2} \Cu^S_{\ast,w_2}$ and $\gr_{\cF_1} \Cu^S_{w_1,\ast}$ denote the graded
complexes associated to the filtrations $\cF_2$ and $\cF_1$, respectively. 
\begin{lemma}\label{griso}
  One has isomorphisms of complexes
  \begin{align}
    \gr_{\cF_2} \Cu^S_{\ast,w_2}&\cong \Cu^S_{\ast,0}\ot M(w_2.0)^\ast \label{grad2}\\
    \gr_{\cF_1} \Cu^S_{w_1,\ast}&\cong \Cu^S_{0,\ast}\ot M(w_1.0) \label{grad1}
  \end{align}
  for $\ast \ge 0$.
\end{lemma}
\noindent{\bf Proof:} For $e\in \sum_{\hght(\beta)\le k}\uqlsp_\beta$ and $u\in \uqg$ one 
obtains in analogy to the
proof of Lemma \ref{filtered}
\begin{align*}
  h_{w_1,w_1';w_2}(u\ot (v_{w_1.0} & \ot e \xi_{-w_2.0}))
             =u e y_{w_1,w_1'}\ot (v_{w_1.0}\ot \xi_{-w_2.0}) \\
             &\in u y_{w_1,w_1'}\ot (v_{w_1.0}\ot e \xi_{-w_2.0}) 
                        +\cF_2^{k{-}1} W(w_1'.0,w_2.0).
\end{align*}
This shows (\ref{grad2}) and (\ref{grad1}) is verified analogously. $\blacksquare$

\begin{proposition}\label{CSexact}
  For any $w_1,w_2\in W^S$ and $n\in \N$ the complexes
  $\Cu^S_{w_1,\ast}$ and $\Cu^S_{\ast,w_2}$ are exact
  in $\Cu^S_{w_1,n}$ and $\Cu^S_{n,w_2}$, respectively.
\end{proposition}
\noindent {\bf Proof:} 
  This follows immediately from Lemma \ref{griso} and the exactness of the complexes
  $\Cu^S_{\ast,0}$ and $\Cu^S_{0,\ast}$.
$\blacksquare$

\vspace{.5cm}

The above proposition is one main step in order to prove that
$\Cu^S_{\ast,\ast}$ defined in (\ref{Ctilde_nm}) together with the maps
\begin{align*}
  \hu_{n,\ast}:\Cu^S_{n,\ast}\rightarrow 
  \Cu^S_{n-1,\ast},\qquad \hu_{n,m}:=
  \sum_{w_2\in W^S,\,l(w_2)=m}\hu_{n,w_2}\\
  \underline{\hb}_{\ast,m}:\Cu^S_{\ast,m}\rightarrow 
    \Cu^S_{\ast,m-1},\qquad \underline{\hb}_{n,m}:=
  (-1)^n \sum_{w_1\in W^S,\,l(w_1)=n}\hu_{w_1,m}\\
\end{align*}
is a double complex. 

\begin{proposition}\label{doubleComplex}
  $(\Cu^S_{\ast,\ast},\hu_{\ast,\ast},
  \underline{\hb}_{\ast,\ast})$ is a double complex, i.e.~for any
  $n,m\in \N$ the relation
  \begin{align}\label{commutes}
    \underline{\hb}_{n-1,m}\circ \hu_{n,m}+
    \hu_{n,m-1}\circ \underline{\hb}_{n,m}=0
  \end{align}
  holds.
\end{proposition}
\noindent {\bf Proof:}
Note first that the claim of the proposition holds for $n=m=1$.
Indeed, recall that the simple reflection $s_s$ corresponding to 
$\alpha_s$ is the only element in $W^S$ of length one and that
$s_s.0=-\alpha_s$. Thus
\begin{align*}
  \Cu^S_{0,1}=\uqg\ot_{\uqls}(M(0)\ot M(-\alpha_s)^\ast)\\
  \Cu^S_{1,0}=\uqg\ot_{\uqls}(M(-\alpha_s)\ot M(0)^\ast)\\
  \Cu^S_{1,1}=\uqg\ot_{\uqls}(M(-\alpha_s)\ot M(-\alpha_s)^\ast)\\
\end{align*}
and up to a sign the maps $\hu_{1,\ast}$ and 
$\underline{\hb}_{\ast,1}$ are 
determined by $y_{s_s,1}=F_s$ and $x_{s_s,1}=E_s$, respectively.
Therefore
\begin{align*}
  (\underline{\hb}_{0,1}\circ \hu_{1,1}+
    \hu_{1,0}\circ \underline{\hb}_{1,1})
    (u\ot (v_{-\alpha_s}\ot \xi_{\alpha_s}))=
    u(F_sE_s-E_sF_s)\ot (v_0\ot\xi_0)=0
\end{align*}
for all $u\in \uqg$.

Now the proof is performed by induction over $n$ and $m$.
Assume that (\ref{commutes}) holds for some $n,m\in \N$.
We will show that this implies (\ref{commutes}) with $n$ replaced
by $n+1$. The induction over $m$ is performed analogously.

Note that (\ref{commutes}) is equivalent to 
\begin{align*}
  0=[y_{w_1,w'_1},x_{w_2,w'_2}]\ot (v_{w'_1.0}\ot \xi_{-w'_2.0})\in 
  \Cu^S_{n-1,m-1}
\end{align*}
for all $w_1,w'_1,w_2,w'_2\in W^S$ such that $l(w_1)=n$, $l(w_2)=m$
and $w_1\rightarrow w'_1$, $w_2\rightarrow w'_2$.
In particular one gets for any $w''_1\in W^S$ such that 
$l(w''_1)=n+1$ the relation
\begin{align*}
  \sum_{\makebox[2cm]{$w_1\in W^S,\, l(w_1)=n \atop 
  w''_1\rightarrow w_1\rightarrow w'_1$}}
  y_{w''_1,w_1}[y_{w_1,w'_1},x_{w_2,w'_2}]
            \ot (v_{w'_1.0}\ot \xi_{-w'_2.0})=0
\end{align*}
and thus using $\hu_{n,\ast}\circ \hu_{n+1,\ast}=0$
one obtains
\begin{align*}
  \sum_{\makebox[2cm]{$w_1\in W^S,\, l(w_1)=n \atop 
    w''_1\rightarrow w_1\rightarrow w'_1$}}
    [y_{w''_1,w_1},x_{w_2,w'_2}]y_{w_1,w'_1}
            \ot (v_{w'_1.0}\ot \xi_{-w'_2.0})=0.
\end{align*}
By the exactness stated in Proposition \ref{CSexact} 
for all $w''_1, w_1, w_2, w'_2\in W^S$, 
$w''_1\rightarrow w_1$, $w_2\rightarrow w'_2$, $l(w_1)=n$, $l(w_2)=m$  
there exist 
elements $u_{w'''_1}\in \uqg$
such that the relation 
\begin{align}\label{rand+1}
  [y_{w''_1,w_1},x_{w_2,w'_2}]\ot (v_{w_1.0}\ot \xi_{-w'_2.0})=
  \sum_{\makebox[0cm]{$w'''_1\in W^S\atop w'''_1\rightarrow w_1$}}
  u_{w'''_1}y_{w'''_1,w_1}\ot (v_{w_1.0}\ot \xi_{-w'_2.0})
\end{align}
holds. By Corollary \ref{notRand} the above relation implies
\begin{align*}
  0=[y_{w''_1,w_1},x_{w_2,w'_2}]\ot (v_{w_1.0}\ot \xi_{-w'_2.0})\in \Cu^S_{n,m-1}.\qquad
\end{align*}
This is relation (\ref{commutes}) with $n$ replaced by $n+1$. $\blacksquare$

\vspace{.5cm}

Using Proposition \ref{doubleComplex} one can now define a double 
complex of covariant $\cB$-bimodules setting
\begin{align*}
  \Om ^{m,n}=
  \bigoplus_{\makebox[0cm]{$w_1,w_2\in W^S\atop l(w_1)=m,\,l(w_2)=n $}}
  \Om (w_1,w_2)
\end{align*}
where for $w_1,w_2\in W^S$ we define
\begin{align*}
  \Om (w_1,w_2):=\{f\in W(w_1.0,w_2.0)^\ast\,|\, \dim (f\uqg)<\infty\}.
\end{align*}
Note that by definition 
$\Om ^{m,n}=\{f\in C_{m,n}\,|\, \dim (f\uqg)<\infty\}$.
Thus the differentials $\hu_{\ast,\ast}$ and 
$\underline{\hb}_{\ast,\ast}$ on $C_{\ast,\ast}$ induce differentials
$\del^{\ast,\ast}$ and $\delb^{\ast,\ast}$ on $\Om ^{\ast,\ast}$,
respectively. Proposition \ref{doubleComplex} implies that
$(\Om ^{\ast,\ast},\del^{\ast,\ast},\delb^{\ast,\ast})$ is a double
complex. Let $(\Om ^\ast, \dif^\ast)$ denote the corresponding 
total complex, i.~e.~ $\Om^n=\bigoplus_{k+l=n}\Om ^{k,l}$ and the
differential $\dif^n:\Om ^{n-1}\rightarrow \Om ^n$ is given by
\begin{align*}
  \dif^n=\sum_{m=0}^n \del^{n-m,m} + \delb^{m,n-m}.
\end{align*}
We are now in a position to formulate the main result of this paper.
\begin{theorem}\label{mainResult}
  There exists an isomorphism $\iota^\ast:\Gdifuw\rightarrow \Om ^\ast$ of complexes
  of covariant $\B$-bimodules.
\end{theorem}
The proof is performed as in Section \ref{delcalc}
up to slight modifications. The details are given for the 
convenience of the reader.

Note first that $(\dif^1:\cB\rightarrow \Om ^1)$ is a covariant
FODC isomorphic to the FODC $(\dif:\cB\rightarrow \Gdifuw[1])$
constructed in \cite{a-heko06}. Indeed, $\Om^1=\Om^{1,0}\oplus
\Om^{0,1}$ and $\Gdifuw[1]=\Gduw[1]\oplus\Gdbuw[1]$ and
the isomorphisms of calculi 
$(\Om^{1,0},\del^{1,0})=(\Gduw[1],\del)$ and 
$(\Om^{0,1},\del^{0,1})=(\Gdbuw[1],\delb)$ have been noted in
subsections \ref{delcalc} and \ref{delbcalc}, respectively.

Next note that the following analogue of Lemma \ref{surlemma} holds.
\begin{lemma}
  The map
    \begin{align*}
      \phi_n:\cB\ot \Om ^n\rightarrow \Om ^{n+1},\qquad 
          b\ot \omega\mapsto b\,\dif^{n+1}\omega
    \end{align*}
  is surjective.
\end{lemma}
\noindent 
{\bf Proof:}

Note that for all $k,l$ the covariant $\cB$-modules $\Om ^{k,l+1}$
and $\Om ^{k+1,l}$ have no irreducible component in common.
Indeed,
\begin{align*}
  (\Xi\circ \Phi)  \Om ^{k,l+1}=
 \bigoplus_{\makebox[0cm]{$w_1,w_2\in W^S\atop l(w_1)=k,\,l(w_2)=l+1$}}
 M(w_1.0)^\ast \ot M(w_2.0)
\end{align*}
lies in the eigenspace corresponding to the eigenvalue
$q^{k-(l+1)}$ of the action of the element $\tau(\omega_s)$ in the
simply connected quantized universal enveloping algebra $\uqgc$ \cite[3.2.9]{b-Joseph}.
Thus to prove surjectivity of $\phi_n$ it is sufficient to show that
the maps
\begin{align*}
  \phi_{k,l}:\cB\ot \Om^{k,l}\rightarrow \Om^{k+1,l},\qquad 
  b\ot \omega\mapsto b \del \omega\\
  \overline{\phi}_{k,l}:\cB\ot \Om^{k,l}\rightarrow \Om^{k,l+1},\qquad 
  b\ot \omega\mapsto b \delb \omega
\end{align*}
are surjective. 

Let $w_1,w_2\in W^S$ such that $\Om (w_1,w_2)\subset \Om^{k+1,l}$
Choose $w'_1\in W^S$ such that $w_1\rightarrow w'_1$. 
Let $\{\xi_i\}$ be a basis of weight vectors of $M(w_2.0)^\ast$.
By the categorical equivalence there exist elements 
$f_i\in \Om(w'_1,w_2)$ such that
\begin{align*}
  f_i(1\ot v_{w'_1.0}\ot \xi_j)=\delta_{ij}.
\end{align*}
Using the fact that $\cB$ separates $\ubar=\uqg/\uqg \uqls^+$ \cite[Proposition~7]{a-heko04} 
and Corollary \ref{uncomparable} one sees that there exists an element
$b_{w_1}\in \cB$  such that
\begin{align*}
  b_{w_1}(y_{w,w'_1(1)})y_{w,w'_1(2)}&=\delta_{w,w_1}
  &&\mbox{for all } w\in W^S,\, w\rightarrow w'_1.
\end{align*}
One obtains
\begin{align*}
  \del^{k,l}(b_{w_1}f_i)(1\ot v_{w.0}\ot \xi_j)=
  (b_{w_1}f_i)(y_{w,w'_1}\ot v_{w.0}\ot \xi_j)=
  \delta_{w,w_1}\delta_{i,j}
\end{align*}
Let $\{v_i\}$ denote a weight basis of $M(w_1.0)$. Acting with
elements of $\uqls$ on $\del^{k,l}(b_{w_1}f_i)$ one obtains
elements $g_{i,j,w_1}$ such that $g_{i,j,w_1}(1\ot v_k\ot \xi_l)=
\delta_{i,k}\delta_{j,l}$ and
\begin{align*}
  g_{i,j,w_1}|_{W(w,w_2)}=0\qquad \mbox{for all } w\in W^S,\,w\neq w_1.
\end{align*}
Thus for the covariant $\cB$-bimodule 
$\Lambda=\im(\phi_{k,l}|_{\cB\ot \Om(w'_1,w_2)})$ the pairing
\begin{align*}
  \Lambda/\cB^+\Lambda\times \bigoplus 
  _{\makebox[0cm]{$w_1\in W^S\atop w_1\rightarrow w'_1$}}M(w_1,w_2)
  \rightarrow \C
\end{align*}
is nondegenerate. By the categorical equivalence one obtains
\begin{align*}
\Lambda=\bigoplus
_{\makebox[0cm]{$w_1\in W^S\atop w_1\rightarrow w'_1$}} \Om (w_1,w_2)
\end{align*}
which proves the surjectivity of $\phi_{k,l}$.
The surjectivity of $\overline{\phi}_{k,l}$ is proved analogously.
$\blacksquare$

\vspace{.5cm}

The remaining steps to identify $(\Om ^\ast, \dif ^\ast)$ with 
$(\Gdifuw,\dif)$ are now straightforward analogues of the 
Lemmata \ref{barwedgeWellDef}, \ref{barwedgeSur}, \ref{ikWellDef} and 
of Proposition \ref{DelCalcIso}. The proofs are omitted.
Define a map
\begin{align}\label{difbarwedgedef}
  \barwedge:\Om ^1\ot_\cB \Om ^n&\rightarrow \Om ^{n+1},\\
  (a\,\dif^0 b)\ot \omega &\mapsto a\,\dif^0 b\barwedge \omega
    :=a\,\dif^{n}(b\omega) - ab\,\dif^{n}\omega.\nonumber
\end{align}
\begin{lemma}\label{difbarwedgeWellDef}
  The map $\barwedge$ is well defined.
\end{lemma}

\begin{lemma}\label{difbarwedgeSur}
  The map $\barwedge:\Om ^1\ot_\cB \Om ^n\rightarrow \Om ^{n+1}$ 
  is surjective.
\end{lemma}
Define a map 
\begin{align*}
  \iota^n:\Gdifuw[n]&\rightarrow \Om ^{n}\\
   a_0\dif a_1\wedge\dif a_2
  \wedge \dots \wedge \dif a_n &\mapsto a_0 \dif^0 a_1\barwedge
  (\dif^0 a_2\barwedge( \dots \barwedge \dif^0 a_n))
\end{align*}  
by repeated use of the map $\barwedge$.
\begin{lemma}\label{difikWellDef}
  The map $\iota^k$ is well defined.
\end{lemma}
\begin{proposition}\label{DifCalcIso}
  The map $\iota^\ast:\Gdifuw[\ast] \rightarrow \Om ^{\ast}$ is an 
  isomorphism of complexes of covariant $\cB$-bimodules.
\end{proposition}

\noindent{\bf Proof of Theorem \ref{mainResult}:} This proof is now
  performed in complete analogy to the proof of Proposition \ref{DelCalcIso}.
  First note that by construction $\iota^\ast:\Gdifuw[\ast] \rightarrow
  \Om ^{\ast}$ is a morphism of complexes. Moreover, by Theorem
  \ref{summaryTheo}(iii), the definition of $\Om^n$, and Proposition
  \ref{C=D} the dimensions of the covariant $\cB$-modules  $\Gdifuw[n]$
  and $\Om ^{k}$ coincide. As $\iota^n$ is a surjection of covariant 
  $\cB$-modules by Lemma \ref{difbarwedgeSur} the categorical 
  equivalence implies that $\iota^n$ is an isomorphism.
  $\blacksquare$

\section{Appendix: Commonly used notation}
  
  Symbols defined in Section \ref{notation1} in order of appearance:
  
  $\gfrak$, $r$, $\hfrak$, $R$, $\pi$, $\alpha_i$, $R^+$, $R^-$, $\nfrak^+$, $\nfrak^-$,
  $(\cdot,\cdot)$, $Q$, $P$, $\alpha_i^\vee$, $\hght$, $\omega_i$, $(a_{ij})$, $P^+$, $Q^+$, $V(\mu)$,
  $\Pi(V(\mu))$, $G$, $S$, $Q_S$, $Q_S^+$, $R_S^\pm$, $P_S$, $P_S^\op$, 
  $\pfrak_S$, $\pfrak_S^\op$, $\lfrak_S$, $L_S$, $P_S^+$, $M(\lambda)$, $W$, $s_\alpha$,
  $W_S$, $W^S$, $l$, $w.\mu$, $\rho$, $w\rightarrow w'$, $w\le w'$.
  
  \medskip
  \begin{tabbing}
     \noindent
       \=Section \ref{uqgcqg}: \qquad\qquad \=\\
       \> $q$ \> element of $\C$, not a root of unity\\
       \> $\uqg$ \> quantized enveloping algebra of $\gfrak$\\
       \> $K_i, K_i^{-1}, E_i, F_i$\> generators of $\uqg$\\
       \> $\kow$, $\kappa$, $\varepsilon$\> coproduct, antipode, and counit of $\uqg$\\
       \> $\ad$\> left adjoint action \\
       \> $\eta$\> algebra isomorphism coalgebra antiautomorphism of $\uqg$\\
       \> $\uqnp$\> subalgebra of $\uqg$ generated by $\{E_i\,|\,1\le i\le r\}$\\
       \> $\uqnm$\> subalgebra of $\uqg$ generated by $\{F_i\,|\,1\le i\le r\}$\\
       \> $U^0$\>  subalgebra of $\uqg$ generated by $\{K_i,K_i^{-1}\,|\,1\le i\le r\}$\\
       \> $G_+$\>  subalgebra of $\uqg$ generated by $\{E_iK_i^{-1}\,|\,1\le i\le r\}$\\
       \> $V(\mu)$\> irreducible left $\uqg$-module of highest weight $\mu\in P^+$\\
       \> $V_\eta$\> $\uqg$-module with action twisted by $\eta$\\
       \> $V(\mu)^\ast$\> right or left $\uqg$-module dual to $V(\mu)$\\
       \> $c_{f,v}^\mu$\> matrix coefficient of $V(\mu)$\\
       \> $C^{V(\mu)}$\> space of matrix coefficients of $V(\mu)$  \\
       \> $\cqg$\> $q$-deformed coordinate ring of $G$\\
       \> \> \\
       \>Section \ref{kebe}: \>\\
       \>$\uqls$ \>subalgebra of $\uqg$ generated by $\{E_i,F_i,K_i^{\pm 1}\,|\,\alpha_i\in S\}$\\
       \>$V_-$ \>subalgebra of $\uqg$ generated by $\{(\ad k)F_i\,|\,k\in \uqls \alpha_i\notin S\}$\\
       \>$V_+$ \>subalgebra of $\uqg$ generated by $\{(\ad k)(E_iK_i^{-1})\,|\,k\in \uqls, \alpha_i\notin S\}$ \\
       \>$\uqlsm$ \> $\uqnm\cap \uqls$\\
       \>$\uqlsp$ \> $G_+\cap\uqls$\\
       \>$\uqps$ \> subalgebra of $\uqg$ generated by $\{E_i,K_i^{\pm 1},F_j\,|\,\alpha_i\in \pi, \alpha_j\in S\}$\\
       \> \> \\
       \>Section \ref{QVnots}: \>\\
       \> $M(\lambda)$\> irreducible left $\uqls$-module of highest weight $\lambda\in P_S^+$ \\
       \> $V^{M(\lambda)}$\> $\uqg\ot_{\uqps} M(\lambda)$ for $\lambda\in P^+_S$\\
       \> $V^\lambda$\> $V^{M(\lambda)}$ for $S=\emptyset$\\
       \> $\rho_S$\> $ \sum_{\alpha\in R_S^+}\alpha/2$\\
       \> \> \\
       \>Section \ref{BGGsection}: \>\\
       \> $C_j^S$\>$\bigoplus_{w\in W^S, l(w)=j}V^{M(w.\mu)}$ for a fixed $\mu\in P^+$ \\
       \> $\varphi^S_j$\> boundary operator of BGG-resolution\\
       \> $f_{w,w'}$\> fixed embedding of Verma modules $V^{w.\mu}\rightarrow V^{w'.\mu}$ if $w\le w'$\\
       \> $s(w_1,w_2)$\> $\pm 1$\\
       \> $h_{w,w'}$\> standard map induced by $s(w,w')f_{w,w'}$\\
       \> $y^\mu_{w,w'}$\> element of $\uqnm$ such that $f_{w,w'}(u\ot v_{w'.\mu})=uy^\mu_{w,w'}\ot v_{w'.\mu}$\\
       \> $y_{w,w'}$\> $s(w,w')y^0_{w,w'}$\\
       \> \> \\
       \>Section \ref{inducedNots}: \>\\
       \> $W^{M(\lambda)}$\> $\uqg\ot_{\uqls} M(\lambda)$ for $\lambda\in P^+_S$\\
       \> $W^\lambda$\> $W^{M(\lambda)}$ for $S=\emptyset$\\
       \> \> \\
       \>Section \ref{functor}: \>\\
       \> $\Phi_\lambda$\> canonical surjection $W^{M(\lambda)}\rightarrow V^{M(\lambda)}$\\
       \> $\Vcal$\> category of finite direct sums of $V^{M(\lambda)}$, $\lambda\in P_S^+$\\
       \> $\Wcal$\> category of finite direct sums of $W^{M(\lambda)}$, $\lambda\in P_S^+$\\
       \> $\fun:\Vcal\rightarrow \Wcal$\> functor defined above Proposition \ref{PresProp} \\
       \> $\uqpsop$\> subalgebra of $\uqg$ generated by $\{E_j,K_i^{\pm 1},F_i\,|\,\alpha_i\in \pi, \alpha_j\in S\}$ \\
       \> \> \\
       \>Section \ref{Homest}: \>\\
       \> $x^\mu_{w,w'}$ \> $\eta(y^\mu_{w,w'})$\\
       \> $\theta_1, \thetab_1, \theta_2, \thetab_2$\> $\uqls$-module homomorphisms defined in and before 
           Proposition \ref{triangProp}\\
       \> $W(\mu,\nu)$\> $\uqg\ot_{\uqls}(M(\mu)\ot M(\nu)^\ast)$\\
       \> $\cF^\ast_1$, $\cF^\ast_2$\> Filtrations of $W(\mu,\nu)$ defined by (\ref{fil1}) and (\ref{fil2})\\
       \> \> \\
       \>Section \ref{Take}: \>\\
       \> $\cA$\> $\cqg$ \\
       \> $\cB$\> $\{b\in \cA\,|\,b_{(1)}b_{(2)}(k)=\vep(k)b\,\mbox{ for all } k\in \uqls \}$\\
       \> $\cB^+$\> $\{b\in \cB\,|\, \vep(b)=0\}$\\
       \> $\Aqr$\> $\cA/\cB^+\cA$\\
       \> $\pair{\cdot}{\cdot}$\> canonical pairing $\uqls\times \Aqr\rightarrow \C$\\
       \> $\AqM$\> category of  finite dimensional left $\Aqr$-comodules \\
       \> $\MK$\> category of finite direct sums of modules of the form $M(\lambda)^\ast$, $\lambda\in P^+_S$\\
       \> $\Xi$\> functor $\AqM\rightarrow \MK$ defined by (\ref{Xi-Def})\\
       \> $P\square_C Q$ \> cotensor product of $P$ and $Q$ over coalgebra $C$ \\
       \> $\ABM$\> category of left $\cA$-covariant left $\cB$-modules\\
       \> $\Phi$\>  functor $\ABM\rightarrow \AqM$ defined by (\ref{Phi-Def})\\
       \> $\Psi$\>  functor $\AqM\rightarrow \ABM$ defined by (\ref{Psi-Def})\\
       \> \> \\
          \>Section \ref{locallyfinite}: \>\\
       \> $\Omega(\lambda)$\>$\{f\in (\WM)^\ast \,|\, \dim (f\uqg)<\infty \}$ \\
       \> $c$\> canonical inclusion $\Om(\lambda) \rightarrow \uqg^\ast$\\
       \> \> \\
       \>Section \ref{q-Diffforms}: \>\\
       \> $(\Gamma_\del,\del), (\Gamma_\delb,\delb)$\> the two irreducible covariant FODC over $\cB$\\
       \> $T_\del, T_\delb$\> quantum tangent space of $\Gamma_\del$ and $\Gamma_\delb$, respectively\\
       \> $(\Gam_\dif,\dif)$\> $(\Gamma_\del\oplus \Gamma_\delb,\del\oplus \delb)$\\
       \> $\Gduw,\Gdbuw,\Gdifuw$\> universal differential calculus of $\Gamma_\del$, $\Gamma_\delb$, 
                                   and $\Gamma_\dif$, respectively\\
       \> \> \\
  \end{tabbing}


\providecommand{\bysame}{\leavevmode\hbox to3em{\hrulefill}\thinspace}
\providecommand{\MR}{\relax\ifhmode\unskip\space\fi MR }
\providecommand{\MRhref}[2]{%
  \href{http://www.ams.org/mathscinet-getitem?mr=#1}{#2}
}
\providecommand{\href}[2]{#2}

\textsc{Istv\'an Heckenberger, Mathematisches Institut, Universit\"at Leipzig,
         Augustusplatz 10-11, 04109 Leipzig, Germany.}
       
 \textit{E-mail address:} \texttt{heckenberger@math.uni-leipzig.de}

\vspace{.5\baselineskip}
                                                
\textsc{Stefan Kolb, Mathematics Department, Virginia Polytechnic
         Institute and State University, Blacksburg, VA 24061 USA}
         
   
\vspace{.5\baselineskip}

\textsc{Present address: School of Mathematics and Maxwell Institute for
Mathematical Sciences, The University of Edinburgh, JCMB, The King's Buildings, Mayfield Road, Edinburgh
         EH9 3JZ, UK}
         
 \textit{E-mail address:} \texttt{stefan.kolb@ed.ac.uk}

\end{document}